
\input amstex
\documentstyle{amsppt}
\pageno=1
\magnification1200
\catcode`\@=11
\def\logo@{}
\catcode`\@=\active
\NoBlackBoxes
\vsize=23.5truecm
\hsize=16.5truecm
\def\d{d\!@!@!@!@!@!{}^{@!@!\text{\rm--}}\!}

\def\inv{^{-1}}

\def\leg{\;\dot{\le}\;}
\def\geg{\;\dot{\ge}\;}

\def\crp{\overline{\Bbb R}_+}

\def\rnp{{\Bbb R}^n_+}

\def\crnp{\overline{\Bbb R}^n_+}

\def\srplus{\Cal S(\overline{\Bbb R}_+)}
\def\srminus{\Cal S(\overline{\Bbb R}_-)}

\def\srnpm{\Cal S(\overline{\Bbb R}^n_\pm)}

\def\ang#1{\langle {#1} \rangle}

\def\Zfrac{\tsize\frac1{\raise 1pt\hbox{$\scriptstyle z$}}}
\def\zfrac{\frac1{\raise 1pt\hbox{$\scriptscriptstyle z$}}}
\def\crpp{\overline {\Bbb R}^2_{++}}
\def\rpp{ {\Bbb R}^2_{++}}
\def\rp{ \Bbb R_+}

\define\tr{\operatorname{tr}}

\define\Tr{\operatorname{Tr}}

\define\srp{\Cal S_+}

\define\srpp{\Cal S_{++}}
\define\stimes{\!\times\!}

\document

\topmatter
\title
{Logarithmic terms in trace expansions of Atiyah-Patodi-Singer problems}
\endtitle
\author{Gerd Grubb}
\endauthor
\affil
{Copenhagen Univ\. Math\. Dept\.,
Universitetsparken 5, DK-2100 Copenhagen, Denmark.
E-mail {\tt grubb\@math.ku.dk}}\endaffil
\rightheadtext{Logarithmic terms}
\abstract
For a Dirac-type operator $D$ on a manifold $X$ with a spectral boundary
condition (defined
by a pseudodifferential projection), the associated heat operator
trace has an expansion in integer and half-integer powers and
log-powers of $t$; the interest
in the expansion coefficients goes back to the work of Atiyah,
Patodi and Singer. In the product case considered by APS, it is known
that all the log-coefficients vanish when $\dim X$ is
odd, whereas the log-coefficients at integer powers vanish when $\dim
X$ is even. We here investigate whether this partial vanishing of
logarithms holds more generally. One type of result, shown for
general $D$ with well-posed boundary conditions, is that a perturbation of
$D$ by a tangential differential operator vanishing to order $k$ on
the boundary leaves the first $k$ log-power terms invariant (and the
non-local power terms of the same degree are only locally
perturbed). Another type of result is that for perturbations of the
APS product case by tangential operators commuting with the
tangential part of $D$, all the logarithmic terms vanish when $\dim X$ is
odd (whereas they can all be expected to be nonzero when $\dim X$ is
even).
The treatment is based on earlier joint work with R\. Seeley and a
recent systematic parameter-dependent pseudodifferential boundary
operator calculus, applied to the resolvent.
\endabstract

\endtopmatter

\subhead 
\footnote{\eightrm  Published in Ann.\ Global An.\ Geom.\ 24, 1--51 (2003)}Introduction \endsubhead

Let  $D$ be a first-order differential operator
of Dirac-type
from $C^\infty (X,E_1)$ to $C^\infty (X,E_2)$ ($E_1$ and
$E_2$ Hermitian $N$-dimensional vector bundles over
a compact $n$-dimensional  $C^\infty $ manifold $X$ with boundary $\partial X=X'$), and let $D_\ge$
be the $L_2$-realization defined by the
boundary condition $\Pi _\ge(u|_{X'})=0$; here $\Pi _\ge$ is the
orthogonal projection onto the nonnegative eigenspace for a certain
selfadjoint operator $A$ over $X'$ entering in $D$.
For $\Delta _B=D_\ge^*D_\ge$ (and likewise for $D_\ge D_\ge^*$),
the following heat trace expansion was
shown in a joint work with Seeley [GS95]: 
$$ 
\Tr (\varphi  e^{-t\Delta _B})\sim
\sum_{-n\le k< 0} a_{ k}t^{\frac{k}2}+ 
\sum_{k\ge 0}\bigl({ a'_{ k}}\log t+{a''_{ k}}\bigr) t^{\frac 
{k}2}\text{ for }t\to 0+.\tag0.1
$$ 
Here $\varphi $ is a smooth morphism in $E_1$; the coefficient $a'_0$
vanishes when $\varphi =1$ near $X'$. The coefficient
$a''_0$ enters in the index of $D_{\ge}$; the geometric
content of the first four $a_k$ (with $k<0$) has been investigated by Dowker,
Gilkey and Kirsten [DGK99], [GK02].

For the case with product structure near $X'$, as studied originally by
Atiyah, Patodi and Singer in [APS75],
the coefficients were described in [GS96] in terms of the
expansion coefficients of zeta and eta functions of $A$. In
particular, it was found that the coefficients $a'_k$
vanish for $k$ even $>0$; moreover, if $n$ is odd, they vanish for
all $k\ge 0$. The remaining coefficients with $k\ne 0$ are nonzero in
general even 
when $\varphi =1$, cf\. Gilkey and Grubb [GG98]. 

We shall here investigate
to what extent this ``partial vanishing of
logarithms'' may hold in non-product cases. Our principal results are:

1) Consider two choices $D_1$ and $D_2$ of $D$, provided with the
same well-posed boundary condition. If they differ by a first-order
tangential differential operator $x_n^lP$ (where $x_n$ is the normal
coordinate), then  the expansions
(0.1) for $D_1$ and $D_2$ have, for $0\le k\le l$, the same log-coefficients 
$a'_k$, and the coefficients $a''_k$ differ only by local terms. In
particular,  the coefficient
$a'_1$ is preserved under perturbations of an operator with product
structure $D^0$ near $X'$ by terms vanishing at $X'$. (Section 3.)

2) If $D$ is a perturbation of the product case $D^0$ (near $X'$)
with $\Pi =\Pi _\ge$, by a
tangential first-order differential operator {\it commuting with} $A$, then
all log-coefficients are zero if $n$ is odd. When $n$ is even,
nontrivial log-terms can in general be expected for both even and odd
$k$. (Section 5.)  

We also derive the related
expansions for resolvent traces and zeta functions, and 
we allow $\varphi $ to be replaced by a differential operator
$F$ (tangential or acting in all variables).
Similar results are shown for the operator families associated with the
eta function.

In preparation for these results, Section 2 gives a review of the
underlying parameter-dependent pseudodifferential boundary operator
calculus, and Section 4 shows the 
structure of the resolvent in the commuting case.

Throughout this paper, $D_1$ and $D_2$ are provided with the same
boundary condition. Perturbations of the boundary condition are
considered e.g\. in [G01$'$] and in [G02].

\subhead
1. Representation formulas \endsubhead

$E_1$
and $E_2$ have 
Hermitian metrics, and $X$ has a smooth volume element, defining 
Hilbert space structures on the sections, $L_2(E_1)$, $L_2(E_2)$. The
restrictions of the $E_i$ 
to the boundary $X'$ are denoted $E'_i$. A neighborhood of
$X'$ in $X$ has the form $X_c=X'\times[0,c[$, and there the $E_i$
are isomorphic to the pull-backs of the $E'_i$. We let $x_n$ denote
the coordinate in $[0,c[$. 
$L_2(E'_i)$) is defined with respect to the volume element $v(x',0)dx'$ on $X'$ induced by the element
$v(x',x_n)dx'dx_n$ on $X$.

When $D$ is a first-order elliptic
differential operator from $C^\infty (E_1)$ to $C^\infty (E_2)$, it
may always be 
written in the following form over $X_c$:
 $$D=\sigma(\partial _{x_n}+A_1),\tag1.1
 $$
where $\sigma $ is a homeomorphism from $E_1|_{X_c}$ to $E_2|_{X_c}$ and $A_1$ for each $x_n$ is an elliptic operator in the $x'$-variable.
We say that $D$ is of {\it product type} when 
$\sigma $ is independent of $x_n$ and is unitary from $E'_1$ to
$E'_2$, and $A_1=A$ independent of $x_n$ and selfadjoint in
$L_2(E_1')$; here the product measure $dx'dx_n$ is used on $X_c$.
We say that $D$ is of {\it non-product type} when 
 $\sigma $ is still independent of $x_n$ and unitary, but the
condition on $A_1$ is relaxed to:
$$ A_1=A+x_nP_1+P_0,\tag1.2
 $$
where $A$ is as above and
the $P_j$  are smooth $x_n$-dependent differential
operators in $x'$ (in short: tangential differential operators) of
order $\le j$. Since 
$P_0=P_0(0)+x_nP'_0(x_n)$ with $P'_0$ of order 0, we may absorb
$x_nP'(0)$ in the term  $x_nP_1$, so {\it we can assume that $P_0$ is
constant in $x_n$ on $X_c$}. 
(In [GS95], [G99], these operators of product type and of non-product type
were said to be ``of Dirac-type''. Some other 
authors restrict that notation to operators that moreover satisfy$$
\sigma ^2=-I,\quad \sigma A=-A\sigma ,\quad D\text{ is selfadjoint on $X$ and $D^2$ principally
scalar},
$$
which we do not assume here.)

Integration by parts shows that the formal adjoint  
$D^*$ equals
$$
D^\ast=(-\partial _{x_n}+A_1')\sigma ^*,\quad A_1'=A+x_n P^*_1+
P'_0, \text{ on }X_c,
$$
where $P'_0=P_0^*-v^{-1}\partial _{x_n}v$. When 
$\partial _{x_n}v(x',0)=0$, $D^*$ may also be written in the form 
$D^\ast=(-\partial _{x_n}+A+x_nP'_1+P_0^*)\sigma ^*$, with
$P'_1-P^*_1$ of order 0.
If $P_0=0$ and $P_1=x_n^{l-1}P_l$ for some $l\ge 1$, so that
$D=\sigma (\partial _{x_n}+A+x_n^lP_l)$, then $D^*$ can be written in
the form $D^\ast=(-\partial _{x_n}+A+x_n^lP'_l)\sigma ^*$ if
$\partial _{x_n}^jv(x',0)=0$ for $1\le j\le l$.

In the product case we often use the notation
$$
D^0=\sigma (\partial _{x_n}+A),\quad {D^0}'=(-\partial
_{x_n}+A)\sigma ^*;\tag1.3
$$these operators have a meaning on $X^0=X'\times\crp$; and ${D^0}'$ is the
formal adjoint of the operator $D^0$ going from $L_2(E^0_1)$ to
$L_2(E^0_2)$, where the $E^0_i$ are 
the liftings 
of the $E'_i$ to $X^0$, and the product measure is used.

By $V_>$, $V_\geq$, $V_<$ or $V_\leq$  we denote the subspaces of $L^2(E'_1)$
spanned by the eigenvectors of $A$ corresponding to eigenvalues which
are $>0$, $\geq0$, $<0$, or $\leq0$. (For precision one can write
$V_>(A)$, etc.) $V_0$
is the nullspace of $A$. 
The corresponding projections are denoted $\Pi_>$, $\Pi _\ge$, etc\. (note
that $\Pi _\ge=\Pi _>+\Pi _0$ and $\Pi _<=I-\Pi _\ge$). They are
pseudodifferential operators 
($\psi $do's) of order 0; $\Pi _0$ has finite rank and is a $\psi $do of order $-\infty $.
We also define$$
A_{\lambda }=(A^2-\lambda )^{\frac12},\text{ for }\lambda \in \Bbb
C\setminus \operatorname{spec}A^2\supset \Bbb C\setminus\crp .\tag1.4
$$
Moreover, we set
$$\aligned
|A|&=(A^2)^\frac12,\quad A'=A+\Pi _0,\quad \text{ so that }|A'|=|A|+\Pi
_0\text{ and}\\
\Pi_>&=\tfrac12\tfrac{|A|+ A}{|A'|}=\tfrac12\tfrac A
{|A'|}+\tfrac12-\tfrac12 \Pi _0, \quad \Pi _\ge = \Pi _>+\Pi _0=
\tfrac12\tfrac A
{|A'|}+\tfrac12+\tfrac12 \Pi _0.
\endaligned\tag1.5
$$

In Sections 4--5 of this paper, we consider the product and 
non-product cases with the
boundary condition$$
\Pi _\ge(A)\gamma _0u=0,\tag1.6
$$
 where $\gamma_0u=u|_{X'}$, defining the realizations $D_{\ge}$ and
$D^0_\ge$, and we denote ${D_\ge}^*D_\ge=\Delta _{B}$,
${D^0_\ge}'D^0_\ge=\Delta ^0_{B}$. However, the more qualitative
results in Section 3  allow the consideration of a general first-order
elliptic operator $D$ with a boundary condition$$
\Pi \gamma _0u=0,\tag1.7
$$
where $\Pi $ is an orthogonal pseudodifferential projection that is {\it
well-posed} with respect to $D$ (cf\. [S69] or [G99]). This means
that when we at each $(x',\xi ')$ in the cotangent sphere bundle of
$X'$ denote by $N^+(x',\xi ')\subset \Bbb C^N$ the space
of boundary values of null-solutions of the model operator (defined
from the principal symbol $d^0$ of $D$),$$
N^+(x',\xi ')=\{\, z(0)\in \Bbb C^N\mid d^0(x',0,\xi ',D_{x_n})z(x_n)=0
,\; z(x_n)\in L_2(\Bbb R_+)^N\,\},
$$
then the
principal symbol $\pi ^0(x',\xi ')$ of $\Pi $ maps $N^+(x',\xi ')$
bijectively onto the range of $\pi ^0(x',\xi ')$ in $\Bbb C^N$. We denote the
realization of $D$ defined by 
(1.7) by $D_{\Pi }$ and again denote ${D_{\Pi }}^*D_{\Pi }=\Delta _B$; it
likewise has a trace expansion (0.1), cf\. [G99].

As shown in [G01$'$], the coefficients
$a_k$ and $a'_k$ with $k\le J-n$
in the trace expansion (0.1) are unaffected by a replacement of $\Pi
$ by a closed range operator $\Pi 
+S$, where $S$ is a pseudodifferential operator of order $\le -J
$ for some $J\ge 1$.

It is explained e.g\. in [GS95], [GS96] how the heat trace expansion (0.1) is
equivalent with the derived resolvent expansion
$$
\Tr (\varphi \partial _ \lambda ^r(\Delta _B-\lambda )^{-1})\sim
\sum_{-n\le k< 0} \tilde a_{ k}(-\lambda )^{-\frac{k}2-r-1}+ 
\sum_{k\ge 0}\bigl({ \tilde a'_{ k}}\log (-\lambda
)+{\tilde a''_{ k}}\bigr)(-\lambda )^{-\frac 
{k}2-r-1},\tag1.8
$$
where $r+1>\frac n2$ and $\lambda \to\infty $ on rays in $\Bbb
C\setminus \Bbb R_+$. The coefficients $\tilde a_k$, $\tilde a'_k$ and $\tilde
a''_k$ are proportional to the coefficients $a_k$, $a'_k$, $a''_k$ in (0.1),
respectively, by universal nonzero proportionality factors (depending
on $r$).
We shall henceforth work with the resolvent. It is well-known that
(0.1) is likewise equivalent with the zeta function expansion$$
\Gamma (s)\Tr (\varphi  \Delta _B^{-s})\sim
\sum_{-n\le k<0} \frac{a_{ k}}{s+\frac{k}2}-\frac{\Tr(\varphi \Pi
_0(\Delta _B))}s+  
\sum_{k\ge 0}\Bigl(-\frac{ a'_{ k}}{(s+\frac
k2)^2}+\frac{a''_{ k}}{s+\frac k2}\Bigr) ,\tag1.9
$$
describing the pole structure of the meromorphic extension of
$\Gamma (s)\Tr (\varphi  \Delta _B^{-s})$ from $\operatorname{Re}s>\frac
n2$ to $s\in\Bbb C$. Here $\Delta _B^{-s}$ is defined by functional
calculus on 
$V_0(\Delta _B)^\perp$ and is taken to be zero on $V_0(\Delta _B)$.
$\Tr (\varphi \Delta _B^{-s})$ is also denoted  
$\zeta (\varphi ,\Delta _B,s)$, the zeta
function.

The coefficients $\tilde a_k,\tilde a'_k,a_k,a'_k$ are locally determined. The
first sum in (0.1), (1.8), (1.9), is sometimes written as a summation
over all $k\ge -n$; we presently use a convention where such local
contributions for $k\ge 0$ are absorbed in the generally nonlocal
coefficients $\tilde a''_k,a''_k$.

There exist several ways of representing the resolvent. 
A direct way is described in [G92]. Another way, introduced in
[GS95] (see also [G99] for the case (1.7)), is to identify the
resolvent as a block in  
the resolvent of an enlarged first-order system, acting in the bundle
$E=E_1\oplus E_2$ over $X$: Let $$
\Cal D=
\pmatrix 0&-{D^*}\\ D&0 
\endpmatrix,\quad
\Cal D_{\Cal B}=
\pmatrix 0&-{D_{\Pi }}^*\\ D_{\Pi }&0 
\endpmatrix  ;\tag1.10
$$
here $\Cal D_{\Cal B}$ is the realization of $\Cal D$ defined by the
boundary condition $$
\Cal B \gamma _0u=0, \text{ with }
\Cal B=\pmatrix \Pi  & \Pi ^\perp\sigma ^* 
\endpmatrix , \quad u=\pmatrix u_1\\u_2 
\endpmatrix .\tag1.11
$$
The operator $\Cal D$ in (1.10) is formally skew-selfadjoint on $X$.
When $\widetilde D$ is an extension of $D$ to an open $n$-dimensional $C^\infty $ manifold
$\widetilde X$ in which $X$ is smoothly imbedded, we define
$\widetilde {\Cal D}$ from $\widetilde D$ as in (1.10) and set, for
$\mu \in\Bbb C\setminus i\Bbb R$, 
$$
\Cal Q_{ \mu } =\pmatrix \mu (\widetilde D^*\widetilde D+\mu
^2)^{-1}&\widetilde D^*(\widetilde D\widetilde D^*+\mu ^2)^{-1} \\ 
-\widetilde D(\widetilde D^*\widetilde D+\mu ^2)^{-1}&\mu (\widetilde D\widetilde D^*+\mu
^2)^{-1}   
\endpmatrix ,\tag1.12
$$
where $(\widetilde D^*\widetilde D+\mu
^2)^{-1}$ (resp\. $(\widetilde D\widetilde D^*+\mu
^2)^{-1}$) is a parametrix of $\widetilde D^*\widetilde D+\mu
^2$ (resp.\ of $\widetilde D\widetilde D^*+\mu
^2$); it can be taken as an inverse when $\widetilde X$ is compact.
Then ${\Cal Q_\mu }$ is a parametrix --- an inverse if
$\widetilde X$ is compact --- of $\widetilde{\Cal D}+\mu $, as is
easily checked.

The operator $\Cal D_{\Cal B}$ is 
skew-selfadjoint as an unbounded
operator in $L_2( E)$, so 
it has a resolvent $\Cal R_\mu =(\Cal D_{\Cal B}+\mu )^{-1}$ for
$\mu \in \Bbb 
C\setminus  i\Bbb 
R$, equal to 
$$
\Cal R_{ \mu } =(\Cal D_{\Cal
B}+\mu )^{-1}=\pmatrix \mu (D_{\Pi}^*D_{\Pi}+\mu
^2)^{-1}&D_{\Pi}^*(D_{\Pi}D_{\Pi}^*+\mu ^2)^{-1} \\ 
-D_{\Pi}(D_{\Pi}^*D_{\Pi}+\mu ^2)^{-1}&\mu (D_{\Pi}D_{\Pi}^*+\mu
^2)^{-1}   
\endpmatrix .\tag1.13
$$
Thus the resolvent $R_\lambda =(D_\Pi^*D_\Pi-\lambda )^{-1}=(\Delta
_B-\lambda )^{-1}$ which we
want to analyze, can be retrieved as 
$$
R_{-\mu ^2} =\mu ^{-1}\pmatrix 1&0 
\endpmatrix \Cal R_\mu \pmatrix 1\\0
\endpmatrix =\mu ^{-1}\Cal R_{\mu ,11},\quad \lambda =-\mu ^2 .\tag1.14
$$
 The resolvent has the structure
$$
R_{\lambda }=Q_{\lambda ,+}+G_\lambda ,\text{ where } Q_{-\mu ^2}=
\mu ^{-1}\pmatrix 1&0  
\endpmatrix \Cal Q _\mu \pmatrix 1\\0
\endpmatrix,  \tag1.15
$$
and $G_\lambda $ is a singular Green operator (more about pseudodifferential
boundary operators in Section 2).

When $D$ is of non-product type and $\Pi =\Pi _\ge$, we consider
along with $\Cal D$ and $\Cal D_{\Cal B}$ the
associated operators of product type$$
\Cal D^0=
\pmatrix 0&-{ D^0}'\\ D^0&0 
\endpmatrix,\quad
\Cal D^0_{\Cal B}=
\pmatrix 0&-{D^0_{\ge}}'\\ D^0_{\ge}&0 
\endpmatrix  ;\tag1.16
$$
here $\Cal D^0$  and $\Cal D^0_{\Cal B}$ act in $E^0=E^0_1\oplus
E^0_2$, and $\Cal D^0_{\Cal B}$ is determined by the same boundary
condition (1.11) as $\Cal D_{\Cal B}$. $\Cal D^0$ extends
to the bundle 
$\widetilde E$ over $\widetilde X^0=X'\times\Bbb R$ obtained by
lifting $E'=E'_1\oplus E'_2$.

In the product situation, we can define the ingredients by functional
calculus from $A$, using the Fourier transform in the $x_n$-variable
only. We can write $Q^0_\lambda =({D^0}'D^0-\lambda
)^{-1}$ as follows: 
$$
Q^0_\lambda =(D_{x_n}^2+A^2-\lambda )\inv=\Cal F^{-1}_{\xi _n\to
x_n}(\xi _n^2+A^2-\lambda )\inv\Cal 
F_{x_n\to \xi _n}.\tag1.17
$$
Moreover, we can describe the boundary operators using the following
notation for the elementary 
Poisson operator $K_{A_\lambda }$, trace operator $T_{A_\lambda }$ of
class 0, and singular Green operator
$G_{A_\lambda }$ of class 0:$$ 
\aligned
K_{A_\lambda }&=
\operatorname{OPK}_n(e^{-x_nA_\lambda
})=\operatorname{OPK}_n\Bigl(\frac1{A_\lambda +i\xi _n}\Bigr),\\
T_{A_\lambda }&=
\operatorname{OPT}_n(e^{-x_nA_\lambda
})=\operatorname{OPT}_n\Bigl(\frac1{A_\lambda -i\xi _n}\Bigr),\\
G_{A_\lambda }&=\operatorname{OPG}_n(e^{-(x_n+y_n)A_\lambda })
=\operatorname{OPG}_n\Bigl(\frac1{(A_\lambda
+i\xi _n)(A_\lambda -i\eta _n)}\Bigr),
\endaligned\tag1.18
$$ 
here we have used both the
symbol-kernel and the symbol notation, with respect to the
$x_n$-coordinate. (One can 
write $(A_\lambda +i\xi _n)^{-1}$ etc\. as fractions in these
formulas, since they are all commuting functions of the selfadjoint
operator $A$.) Explicitly, for $v\in C^{\infty }(E'_1)$, $u\in C^\infty
(E^0_1)$ with compact support in $X_c$,
$$
[K_{A_\lambda }v](x)=e^{-x_nA_\lambda }v(x'),\quad [T_{A_\lambda
}u](x')=\int_0^\infty e^{-x_nA_\lambda }u(x',x_n)\,dx_n,\tag1.19 
$$
in the sense of functional calculus, and
$$
G_{A_\lambda }u=K_{A_\lambda }T_{A_\lambda }u .\tag1.20
$$

The use of formulas based on functional calculus will be pursued in Section 4.

\example{Remark 1.2} For a resolvent $(T-\lambda )^{-1}$, one has that 
$$
\partial _\lambda ^r(T-\lambda )^{-1}=r!(T-\lambda )^{-r-1},\tag1.21
$$
so it makes no difference whether one refers to powers or to $\lambda
$-derivatives when describing trace expansions for iterated resolvents.
However, in [GS95], the variable $-\lambda $ was replaced by $\mu ^2$
and the primary results were expressed for $\mu $-derivatives, since
this looked less complicated than a description of the many terms resulting
from raising the resolvent to a power. There
were not given many details on how one gets back to the desired
expansions of $\lambda $-derivatives. In fact, a direct consideration
of powers (still departing from a reformulation in the variable $\mu $) would
have been more 
adequate; this road was followed in subsequent treatments (partly in [G99],
fully in [G01]), and is also followed below. The difference lies in
the fact that, with 
notation as (1.14), the power formula $R_\lambda ^2=\mu ^{-2}\Cal
R_{\mu ,11}^2$ shows a decrease in the order (for fixed $\mu $)
which is harder to see 
from the differentiation formula 
$$
\partial _\lambda R_\lambda =\partial _\mu R_{-\mu ^2}\partial
_\lambda \mu =\partial _\mu (\mu ^{-1}\Cal R_{\mu ,11})\tfrac{-1}2
\mu ^{-1}=c_1\mu ^{-3}\Cal R_{\mu ,11}+c_2\mu ^{-1}\partial _\mu \Cal
R_{\mu ,11}.\tag1.22
$$
\endexample

\subhead
2. The parameter-dependent symbol calculus
\endsubhead
 
Let us briefly recall the symbol spaces for pseudodifferential
boundary operators \linebreak($\psi $dbo's) introduced in [G01].
It is simplest to explain for operators of class 0, which is
essentially all we need here (we refer to [G01] for the full
calculus). 

One considers systems of $\mu $-dependent operators$$
\Cal A(\mu )=\pmatrix P(\mu )_++G(\mu )& K(\mu )\\\quad\\ T(\mu )& Q(\mu )
\endpmatrix \:\matrix C^\infty (\crnp)^N &\quad& C^\infty
(\crnp)^{N'}\\
\times&\to&\times\\
C^\infty (\Bbb R^{n-1})^M &\quad & C^\infty (\Bbb R^{n-1})^{M'}
\endmatrix,\tag2.1
$$
which for each fixed $\mu $ belong to the calculus of Boutet de Monvel
[BM71]:
$P(\mu )$ is a $\psi $do on $\Bbb R^{n}$ satisfying the transmission
condition at $x_n=0$, $G(\mu )$ is a singular Green operator
(s.g.o.), $T(\mu )$ is a trace 
operator, $K(\mu )$ is a Poisson operator and $Q(\mu )$ is a $\psi
$do on $\Bbb R^{n-1}$.
For the reader who is not familiar with this calculus, we refer to
e.g\. [G96, Ch\. 1] or [G01, Sect\. 1].

The starting point in the present parameter-dependent case is the
$\psi $do symbol spaces 
from [GS95], based on $x'\in {\Bbb R}^{n-1}$ and now allowed to take
values in  Banach spaces $B$ such as $L_p(\rp)$, $L_p(\rpp)$
(we write $\rpp=\rp\stimes\rp$). 
Moreover, we now take
powers of $|(\xi ',\mu )|$ into the 
definition. 
To smooth out the behavior of $|(\xi ',\mu )|$ near 0, it is
convenient to replace it by $[(\xi ',\mu )]$; 
here $[x]$ denotes a $C^\infty $ function of $x\in{\Bbb
R}^N $ satisfying $[x]=|x|$ for $|x |\ge 1$, $[x]\in [\tfrac12,1]$ for 
$|x|\le 1$. We also use the notation $\ang x=(1+|x|^2)^\frac12$. 
Note that $[(\xi ',{1/z})] =|(\xi ',{1/z})|$ for $|z|\le 1$, and that
$$
|(\xi ',{1/z})|=|(\xi ',\mu )| =|z|^{-1}
\ang{z\xi '},\text{ when }\mu ={1/z}
.\tag2.2
$$
$[(\xi ',\mu )]$ is more briefly written $[\xi ',\mu ]$; it will in the
following often be denoted $\kappa $
(as in [G96]), so from now on,
$$ 
\kappa =[\xi ',{1/z}]=[\xi ',\mu ],\text{ with }\mu ={1/z}.\tag2.3
$$ 
We denote $\{0,1,2,\dots\}=\Bbb N$.
\proclaim{
 Definition 2.1}  Let $m\in {\Bbb R}$, $d$ and $s \in {\Bbb Z}$.
Then $S^{m,0,0} ({\Bbb R}^{n-1}\stimes {\Bbb R}^{n-1},\Gamma 
,B)$ consists of the $C^\infty $ functions $p(x',\xi ',\mu )$ valued
in $B$ which
 satisfy, with $1/\mu  =z$,
$$\multline 
\partial 
_{|z|}^jp(\cdot,\cdot, {1/z} )\in S^{m+ j}({\Bbb R}^{n-1}
\stimes {\Bbb R}^{n-1},B)\text{ for }{1/z} \in\Gamma ,\\
 \text{with uniform estimates for }|z|\le 1,
{1/z}  \text{ in closed
subsectors of }\Gamma ,
\endmultline
\tag 2.4
$$
for all $j\in\Bbb N$.
Moreover, we define 
$$
S^{m,d,s } ({\Bbb R}^{n-1}\stimes {\Bbb R}^{n-1},\Gamma ,B)
=\mu ^d [\xi ',\mu ]^s
S^{m,0,0}({\Bbb R}^{n-1}\stimes {\Bbb R}^{n-1},\Gamma,B
).
\tag2.5
$$ 
\endproclaim 

The indication $ ({\Bbb R}^{n-1}\stimes {\Bbb R}^{n-1},\Gamma
,B)$ is often abbreviated to $ (\Gamma ,B)$, or just $(\Gamma )$ if
$B=\Bbb C$.
Keeping the identification of $\mu $ with ${1/z}$ in mind, we shall
also say that $p(x',\xi ',{1/z})$ lies in $S^{m,d,s}(\Gamma ,B)$.

We leave the requirement (from [GS95]) of being holomorphic in 
$\mu \in \Gamma ^\circ$ out of the definition
since $\kappa =[\xi ',\mu ] $ is not so (one could instead work with a
variant of $\kappa $ that is holomorphic on suitable sectors, as in
[G96], (A.2$''$)--(A.2$'''$)). The symbol is just assumed to be $C^\infty
$ in $\mu  
\in\Gamma $ considered as a subset of ${\Bbb R}^2$.
Accordingly, we write $\partial _{|z|}$ (instead of $\partial _z$),
since it is a control of the radial derivative that is needed
(uniformly when the argument of $z$ runs in a compact interval), and 
$|z|$ in the following enters as a real parameter.

To define symbol-kernels for the boundary operators, we use the cases
$B=L_\infty (\rp)$ and $B=L_\infty (\rpp)$, with variables $x_n$ or
$u_n$, resp.\ $(x_n,y_n)$ or 
$(u_n,v_n)$; these variables will then be mentioned in the detailed
description of the function.

We denote by $r^\pm$ the restriction from distributions on
$\{x_n\in\Bbb R\}$ to distributions on
$\{x_n\gtrless 0\}$, and by $e^\pm$ 
the extension from functions on $\{x_n\gtrless 0\}$ to functions on
$\{x_n\in\Bbb R\}$  
by assigning zero values on $\{x_n\lessgtr 0\}$, respectively. With
$\Cal S(\Bbb R^n)$ denoting the Schwartz space, we denote
$r^\pm\Cal S(\Bbb R^n)=\srnpm$, $\Cal S(\overline{\Bbb R}_\pm)=\Cal S_\pm$. We denote
$r^+_{x_n}r^+_{y_n}\Cal S(\Bbb R^2)=\Cal S_{++}$.

Let us also recall the notation developed from [BM71]:
For $n=1$, the Fourier transformed spaces are denoted  
$\Cal H^+=\Cal F e^+ \srplus$, $\Cal H^-_{-1}=\Cal F e^-
\srminus$; they consist of $C^\infty $ functions that extend
holomorphically to 
$t\in \Bbb C_-$ resp\. $\Bbb C_+$ ($\Bbb C_\pm=\{z\in\Bbb C\mid
\operatorname{Im}z\gtrless 0\}$) and are $O(t^{-1})$ there. Adding to
$\Cal H^-_{-1}$ the space 
$\Bbb C[t]$ of polynomials in $t$, we get the space $\Cal H^-$. We
denote
$$
\Cal H=\Cal
H^+\dot+\Cal H^-,\text{ with projections }h^\pm\:\Cal H\to \Cal
H^\pm.\tag2.6
$$
The space $\Cal F_{x_n\to\xi _n}\overline{\Cal F}_{y_n\to\eta
_n}e^+_{x_n}e^+_{y_n}\Cal S_{++}$ identifies with $\Cal
H^+\hat\otimes\Cal H^-_{-1}$. (Our
notation for the Fourier transform and the conjugate Fourier
transform is:$$
(\Cal F f)(\xi )=\hat f(\xi )=\int_{\Bbb R^n}e^{-ix\cdot\xi
}f(x)\,dx,\quad (\overline{\Cal F} f)(\xi )=\int_{\Bbb R^n}e^{+ix\cdot\xi
}f(x)\,dx,\tag2.7
$$
so that $\Cal F^{-1}=(2\pi )^{-n}\overline{\Cal F}$; they are
sometimes just applied in the $x'$-variable or the $x_n$-variable
alone.)

The appropriate definition of parameter-dependent symbol-kernels for
boundary operators involves a scaling in the $x_n$-variable:

\proclaim{
Definition 2.2}
Let $m  \in{\Bbb R}$, $d$ and $s\in {\Bbb Z}$.

{\rm (i)} The space
${\Cal S}^{m,d,s } ({\Bbb R}^{n-1} \stimes{\Bbb R}^{n-1},\Gamma, \srp
)$ (briefly denoted ${\Cal S}^{m,d,s}  (\Gamma ,\srp)$)  
consists of the  complex
functions $ \tilde f({x'},x_n,{\xi '},\mu )$ in $C^\infty 
({\Bbb R}^{n-1} \stimes\crp\stimes{\Bbb R}^{n-1}\stimes\Gamma )$ satisfying, for all
$l,l'\in{\Bbb N}$,
$$
\ang{z\xi '}^{l-l'}u_n^l\partial _{u_n}^{l'} \tilde f({x'},|z|u_n,{\xi
'},{1/z} 
) \in S^{m,d,s +1}  ({\Bbb R}^{n-1} \stimes{\Bbb R}^{n-1},
\Gamma , L_{\infty ,u_n}(\rp))
\tag2.8$$
(equivalently, $u_n^l\partial _{u_n}^{l'} \tilde f({x'},|z|u_n,{\xi
'},{1/z} 
)$ belongs to $\in S^{m,d+l-l',s+1-l+l' }  (\Gamma , L_{\infty ,u_n}(\rp))$).

{\rm (ii)}
The space
${\Cal S}^{m,d,s} ({\Bbb R}^{n-1}\stimes{\Bbb R}^{n-1},\Gamma ,\srpp)$
(briefly denoted ${\Cal S}^{m,d,s} (\Gamma ,\srpp)$)  
consists of the  complex
functions $ \tilde f({x'},x_n,y_n,{\xi '},\mu )$ in $C^\infty 
({\Bbb R}^{n-1} \stimes\crpp\stimes{\Bbb R}^{n-1}\stimes\Gamma )$
satisfying, for all 
$l,l',k,k'\in{\Bbb N}$,
$$\multline
\ang{z\xi '}^{l-l'+k-k'}u_n^l\partial _{u_n}^{l'} v_n^{k}\partial
_{v_n}^{k'}\tilde f({x'},|z|u_n,|z|v_n,{\xi '},{1/z}
)\\
\in  S^{m,d,s+2 }  ({\Bbb R}^{n-1} \stimes{\Bbb
R}^{n-1},\Gamma , L_{\infty ,u_n,v_n}(\rpp)). 
\endmultline\tag2.9$$
\endproclaim 

In details, the statement in (2.8) means that for all $j$,
$$
\| \partial 
_{|z|}^j(z^{d}\kappa ^{-s-1}\ang{z\xi '}^{l-l'}u_n^l\partial _{u_n}^{l'}\tilde
f({x'},|z|u_n,{\xi '},{1/z})) \|_{L_{\infty ,u_n}}
 \leg\ang{\xi 
'}^{m+j} ,\tag2.10
$$ 
with similar estimates for derivatives $\partial _{\xi '}^\alpha
\partial _{x'}^\beta $ with $m$ replaced by $m-|\alpha |$. There is a
related explanation of (2.9). We here use $\leg$ to indicate
``$\le$ a constant times''; also $\geg$ will be used, and $\dot=$ indicates
that both $\leg$ and $\geg$ hold.

The third upper index $s$ is included to keep track of factors $\kappa
=[\xi ',{1/z}]$ in
a manageable way. When $s=0$, we may lave it out of the notation,
consistently with [GS95]:$$
S^{m,d,0}=S^{m,d}.\tag 2.11
$$ For
the trace formulas later on, it is important to know that we always
have inclusions (that follow from [GS95, Lemma 1.13]):
$$ 
\aligned
S^{m,d,s}&\subset S^{m+s,d,0}\cap
S^{m,d+s,0}\text{ if }s\le 0,\\ 
S^{m,d,s}&\subset S^{m+s,d,0}+
S^{m,d+s,0}\text{ if }s\ge 0,
\endaligned\tag2.12$$
for $S$- as well as for $\Cal S$-spaces. We denote $$
\bigcap_{m\in\Bbb
R}S^{m,d,s}=S^{-\infty ,d,s}, \quad \bigcup_{m\in\Bbb
R}S^{m,d,s}=S^{\infty ,d,s},\quad\text{ etc.;}
\tag2.13$$  
observe that by (2.12), $S^{m ,d,-\infty }=S^{-\infty  ,-\infty
,-\infty }$.      

The following rule follows from the definition (proof details are
given in [G01, Lemma 2.10]):

\proclaim{Lemma 2.3}

{\rm (i)}
When $\tilde f\in \Cal S^{m,d,s}  (\Gamma,\srp )$, then
$x_n^j\partial _{x_n}^{j'}\tilde f\in \Cal S ^{m,d,s-j+j'}
(\Gamma,\srp )$ for $j,j'\in\Bbb N$.

{\rm (ii)} When $\tilde f\in \Cal S^{m,d,s}(\Gamma,\srpp )$, then
$x_n^i\partial _{x_n}^{i'}y_n^{j}\partial _{y_n}^{j'}\tilde f
\in \Cal
S^{m,d,s-i+i'-j+j'}(\Gamma,\srpp )$ for $i$, $i'$, $j$, $j'\in\Bbb N$.
\endproclaim 

In the applications to trace formulas, the symbols moreover have to be
holomorphic in $\mu $ for $\mu \in\Gamma ^{\circ}$ with $|(\xi ',\mu
)|\ge \varepsilon $ (some $\varepsilon >0$); we call such symbols
holomorphic in $\mu $, and
this property is preserved in compositions.

\proclaim{Definition 2.4}

$1^\circ$ The functions in ${\Cal S}^{m,d,s }(\Gamma ,\srp)$ 
are the {\bf Poisson symbol-kernels} and {\bf trace symbol-kernels of
class 0}, 
of degree $m+d+s$, in the parametrized calculus.

$2^\circ$ The functions in ${\Cal S}^{m,d,s }(\Gamma ,\srpp)$
are the
{\bf singular Green symbol-kernels of
class $0$}
and degree $m+d+s$ in the parametrized calculus. 

\endproclaim 

Operators are defined from these symbol-kernels as follows:
$$\aligned
T&=\operatorname{OPT}(\tilde f)\:u(x)\mapsto \int_{\Bbb
R^{2(n-1)}}\int_0^\infty e^{i(x'-y')\cdot
\xi '} \tilde 
f(x',x_n,\xi ',\mu )u(y',x_n)\,dx_ndy'\d\xi ',\\
K&=\operatorname{OPK}(\tilde f)\:v(x')\mapsto \int_{\Bbb
R^{2(n-1)}}e^{i(x'-y')\cdot
\xi '} \tilde 
f(x',x_n,\xi ',\mu )v(y')\,dy'\d\xi ',\\
G&=\operatorname{OPG}(\tilde f)\:u(x)\mapsto \int_{\Bbb
R^{2(n-1)}}\int_0^\infty e^{i(x'-y')\cdot
\xi '} \tilde 
f(x',x_n,y_n,\xi ',\mu )u(y)\,dy\d\xi ';
\endaligned\tag2.14$$
note that the usual $\psi $do definition is used with respect to the
$x'$-variable. Here
$\d \xi '$ stands for $(2\pi )^{1-n}d\xi '$. When these
definitions are applied with respect to the $x_n$-variable alone, we
write $\operatorname{OPT}_n$, $\operatorname{OPK}_n$, $\operatorname{OPG}_n$.

Let us also mention the definition of the associated symbols, where a
Fourier transformation has been performed in the $x_n$-variable. 

\proclaim{Definition 2.5}

{\rm (i)} By Fourier transformation in $x_n$, $e^+{\Cal S}^{m,d,s
}(\Gamma ,\srp)$ is carried over to the 
space \linebreak${\Cal S}^{m,d,s }(\Gamma ,\Cal H^+)$  of {\bf Poisson symbols}
of degree $m+s+d$, and conjugate Fourier transformation in $x_n$
gives the space ${\Cal S}^{m,d,s 
}(\Gamma ,{\Cal H}^-_{-1})$  of {\bf trace symbols of class }0
and degree $m+s+d$.

{\rm (ii)} By Fourier transformation in $x_n$ and conjugate Fourier
transformation in 
$y_n$, \linebreak$e^+_{x_n}e^+_{y_n}{\Cal S}^{m,d,s }(\Gamma ,\srpp)$ is
carried over to the space ${\Cal S}^{m,d,s
}(\Gamma ,{\Cal H}^+\hat\otimes {\Cal H}^-_{-1})$ of 
{\bf singular Green symbols of
class $0$} and degree $m+d+s$. 
\endproclaim 

The $\mu $-dependent 
$\psi $do's $P$ given on $\Bbb R^n$ should in addition to the conditions in
Definition 2.1 satisfy an appropriate transmission condition at
$x_n=0$, which assures that the truncated $\psi $do $P_+=r^+Pe^+$
enters in the calculus in a good way. The general condition is
explained in [G01, Sect\. 6], where 
the class of such symbols, of degree $s$, is denoted 
$S^{0,0,s}_{\operatorname{ut}}(\Bbb R^n\stimes\Bbb R^n,\Gamma )$
(here ``ut'' stands for ``uniform transmission condition'', cf\. also
[G96]).

When $\tilde f_{m-j}\in S^{m-j,d,s } (\Gamma , \Cal S_+)$ for
$j\in\Bbb N$ and $\tilde f\in S^{m,d,s } (\Gamma ,\Cal S_+)$, we say that
$\tilde f\sim\sum_{j}\tilde f_{m-j}$ in 
$S^{m,d,s }(\Gamma ,\Cal S_+) $ if
$\tilde f-\sum_{j<J}\tilde f_{m-j}\in S^{m-J,d,s }(\Gamma ,\Cal S_+)$ for any
$J\in\Bbb N$. 
For any given sequence $\tilde f_{m-j}\in S^{m-j,d,s } (\Gamma ,\Cal
S_+)$, one can construct an $\tilde f$ such that $\tilde
f\sim\sum_j\tilde f_{m-j}$ in $S^{m,d,s}(\Gamma ,\Cal S_+)$.
(Similar statements hold with $\Cal S_{++}$.) 

Of particular interest are the subspaces of the above symbol-kernel
spaces consisting 
of the functions $\tilde f\in S^{m,d,s}$ that are asymptotic series of terms
$\tilde f_{m-j}\in S^{m-j,d,s}$ with a specific
quasi-homogeneity in $(x_n, 
\xi ',\mu )$ or $(x_n,y_n,\xi ',\mu )$, the corresponding Fourier transformed
terms being ordinarily homogeneous in $(\xi _n,\xi ',\mu )$ resp.\ 
$(\xi _n,\eta _n,\xi ',\mu )$ for $|\xi '|\ge c>0$, of 
degree $m-j+d+s$. Such
symbol-kernels and symbols are called (weakly) polyhomogeneous.

The explanation for the $+1$ resp.\ $+2$ in the third upper index in (2.8)
resp.\ (2.9)  is, that with this choice,
$m+d+s$ is consistent with the top degree of homogeneity in the
Fourier transformed 
situation for polyhomogeneous symbols, where the scalings in $x_n$
and $y_n$ lead to  shifts in the 
indices. Further details in [G01].

The label {\it strongly polyhomogeneous} is reserved for those
symbol-kernels and symbols for which
the terms have the homogeneity property on the larger set where
$|(\xi ', \mu )|\ge c$ ($\xi '\in\Bbb R^{n-1}$, $\mu \in \Gamma \cup
\{0\}$), and standard estimates for symbols-kernels in one 
more cotangent variable hold when the extra variable is identified with $|\mu |$ on each ray in
$\Gamma $. Such symbol-kernels and symbols form a subset of the
weakly polyhomogeneous 
symbol-kernel and symbol spaces 
by [GS95, Th\. 1.16] and [G01, Th\. 3.2]:


\proclaim{Theorem 2.6} (On strongly polyhomogeneous symbol-kernels
and symbols.) 

{\rm (i)} When $p$ is a standard polyhomogeneous $\psi $do symbol of
degree $m$ with respect
to $n$ variables (with global estimates), then the symbol obtained by
fixing $x_{n}$ and replacing $\xi _{n}$ by $\mu \in \rp$ is in
$S^{0,0,m}(\Bbb R^{n-1}\stimes\Bbb R^{n-1},\rp)$. 

{\rm (ii)} When $p$ is as in {\rm (i)} with $n$ replaced by $n+1$ and
satisfies the uniform transmission  
condition at $x_n=0$, with respect to the $x_n$-variable, then the
symbol obtained by fixing $x_{n+1}$ and 
replacing $\xi _{n+1}$ by $\mu \in \rp$ is in 
$S^{0,0,m}_{\operatorname{ut}}(\Bbb R^n\stimes\Bbb R^n,\rp)$.

{\rm (iii)} When $\tilde t$, $\tilde k$ or $\tilde g$ is
a polyhomogeneous trace, Poisson or 
singular Green symbol-kernel with respect to $n+1$ variables
$(x_1,\dots,x_{n+1})$ in 
the standard $\psi $dbo calculus of degree $m$ (with global estimates), $x_n$
denoting the normal variable, then the symbol-kernels obtained by
fixing $x_{n+1}$ and replacing $\xi _{n+1}$ by $\mu \in \rp$ are in
$\Cal S^{0,0,m}(\Bbb R^{n-1}\stimes\Bbb R^{n-1},\rp,\Cal S_+)$
resp\. $\Cal S^{0,0,m}(\Bbb R^{n-1}\stimes\Bbb R^{n-1},\rp,\Cal S_{++})$.

\endproclaim 

\demo{Proof} (i) follows by applying [GS95, Th\.
1.16] to the symbol $[\xi ]^{-m}p$ of degree $0$.

As for (ii), the uniform transmission condition is defined for
parameter-dependent symbols in [G01, Sect\. 6]
precisely so that it holds in this situation.

(iii) is shown in [G01, Th\. 3.2].
\qed
\enddemo 

The class of strongly polyhomogeneous $\psi $do symbols satisfing the uniform
transmission condition as in this theorem, on each ray in $\Gamma $,
is denoted  
 $\Cal S^{0,0,m}_{\operatorname{sphg,ut}}(\Gamma )$.

For a simple example of a Poisson operator as in the theorem, see
[G01, Ex\. 3.4].

The symbol-kernel spaces can of course also be defined
for $x'$ runnning in open subsets $U'$ of $\Bbb R^{n-1}$; then
$\Bbb R^{n-1}\stimes\Bbb R^{n-1}$ is replaced by $U'\stimes \Bbb
R^{n-1}$ in the formulas in Definitions 2.1 and 2.2. Likewise, $\Bbb
R^n\stimes\Bbb R^n$ can be replaced by $U\stimes\Bbb R^n$ in the
definition of $\psi $do symbols. The cotangent variable $\xi '$ need
only run in a conical subset of $\Bbb R^{n-1}$. The operators and symbols behave in
a standard way under coordinate transformations (one just has to keep
check of the uniformity in $z$ of the relevant estimates); we shall not
give any details here but just mention that this allows the definition
of operators acting in vector bundles over manifolds, by use of local
coordinates and local trivializations.

The following composition rules are proved in [G01, Ths\. 6.7--6.9]
(recalled here primarily for operators of class zero):

\proclaim{ Theorem 2.7} 
Let $P$ ($\psi $do on $\Bbb R^n$), $G$ , $T$ and  $K$ (class $0$ singular
Green, trace resp\. Poisson operator for $\rnp$),  and $Q$ ($\psi $do
on $\Bbb R^{n-1}$) be parameter-dependent with
sym\-bol(-ker\-nels) $p, \tilde g, \tilde t, \tilde k, q$ satisfying (for
some $m,d,s\in\Bbb Z$): 
$$
\aligned
p(x,\xi ,\mu )&\in S^{0,0,s}_{\operatorname{ut}}(\Bbb R^n\stimes\Bbb R^n,\Gamma ),\\
\tilde g(x',x_n,y_n,\xi ',\mu )&\in 
{\Cal S}^{m,d,s}  (\Bbb R^{n-1}\stimes\Bbb R^{n-1},\Gamma,\Cal S_{++}
),\\
\tilde t(x',x_n,\xi ',\mu ),
\tilde k(x',x_n,\xi ',\mu )&\in 
{\Cal S}^{m,d,s} (\Bbb R^{n-1}\stimes\Bbb R^{n-1},\Gamma,\Cal S_+
),\\
q(x',\xi ',\mu )&\in S^{m,d,s}(\Bbb R^{n-1}\stimes\Bbb R^{n-1},\Gamma ),
\endaligned\tag2.15
$$
and let $P'$, $G'$, $T'$, $K'$ and $Q'$ be given
similarly with $m$, $d$, $s$ replaced by $m'$, $d'$,
and $s'$. Define
$$
m''=m+m',\quad
d''=d+d',\quad
 s ''=s +s '.
\tag2.16$$
Assume that $s$ resp\. $s'$ is $\le 0$ in the formulas where $P$
resp\. $P'$ enter. Then (omitting the indication $\Bbb
R^{n-1}\stimes\Bbb R^{n-1}$): 

{\rm (i) } $TP'_+ $ (trace operator) has symbol-kernel  in $\Cal
S^{m,d,s''}(\Gamma , \Cal S_{+})$,

{\rm (ii) }$\gamma _0P'_+ $ (trace op\.) has symbol-kernel  in $\Cal
S^{0,0,s'}(\Gamma , 
\Cal S_{+})$,

{\rm (iii) }$P_+K' $ (Poisson op\.) has symbol-kernel  in $\Cal
S^{m',d',s''}(\Gamma , \Cal S_{+})$,

 {\rm (iv) }$P_+G' $ (s.g.o\.) has symbol-kernel  in $\Cal
S^{m',d',s''}(\Gamma , \Cal S_{++})$,
 
{\rm (v) }$GP'_+ $ (s.g.o\.) has symbol-kernel  in $\Cal
S^{m,d,s''}(\Gamma , 
\Cal S_{++})$,
 
{\rm (vi) }$GG' $ (s.g.o\.) has symbol-kernel  in $\Cal
S^{m'',d'',s''+1}(\Gamma , \Cal S_{++})$,
 
{\rm (vii) }$KT' $ ($\psi $do) has symbol-kernel  in $\Cal
S^{m'',d'',s''}(\Gamma , \Cal S_{++})$,

{\rm (viii) }$TG' $ (trace op\.) and $GK'$ (Poisson op\.) have symbol-kernels  in $\Cal
S^{m'',d'',s''+1}(\Gamma , 
\Cal S_{+})$,
 
{\rm (ix) }$\gamma _0G' $ (trace op\.) has symbol-kernel  in $\Cal
S^{m',d',s'+1}(\Gamma , 
\Cal S_{+})$,
 
{\rm (x) }$TK' $ ($\psi $do) has symbol  in $S^{m'',d'',s''+1}(\Gamma )$,
 
{\rm (xi) }$\gamma _0K' $ ($\psi $do) has symbol  in $
S^{m,d,s''+1}(\Gamma )$,

{\rm (xii) }$QT'$ (trace op\.) and $KQ' $ (Poisson op\.) have symbol-kernels  in $\Cal
S^{m'',d'',s''}(\Gamma , 
\Cal S_{+})$,

{\rm (xiii) }$QQ' $ ($\psi $do) has symbol  in $
S^{m'',d'',s''}(\Gamma )$,

{\rm (xiv) }$ P_+P'_+=(PP')_+-G^+(P)G^-(P'), $ where 
$ PP'$ ($\psi $do) has symbol  in \linebreak$
S^{0,0,s''}_{\operatorname{ut}}(\Bbb R^n\stimes\Bbb R^n,\Gamma )$,
$G^+(P)$ (s.g.o\.) has symbol-kernel  in $\Cal
S^{0,0,s-1}(\Gamma , \Cal S_{++})$, and
$G^-(P')$ (s.g.o\.) has symbol-kernel  in $\Cal
S^{0,0,s'-1}(\Gamma , \Cal S_{++})$.
\endproclaim

Observe the general principle that the third upper index is lifted
by 1 in the cases 
where the composition involves an integration in $x_n$.

When $\tilde g$ is a singular Green symbol-kernel, we define the {\it normal 
trace} by$$
(\operatorname{tr}_n\tilde g)(x',\xi ',\mu )=\int_0^\infty \tilde g(x',x_n,x_n,{\xi '},\mu
)\,dx_n\tag 2.17$$ 
This a $\psi $do symbol in the calculus: 

\proclaim{Proposition 2.8} When 
$\tilde g(x',x_n,y_n,{\xi '},\mu )
\in \Cal
S^{m,d,s-1 }(\Bbb R^{n-1}\stimes\Bbb
R^{n-1},\Gamma ,\srpp)$,
then the normal trace of $\tilde g$ is a $\psi $do
symbol in $S^{m,d,s}(\Bbb R^{n-1}\stimes\Bbb
R^{n-1},\Gamma ,\Bbb C.)$
\endproclaim 

When $G=\operatorname{OPG}(\tilde g)$, we denote the $\psi $do with
symbol $\tr_n\tilde g$ by $\tr_nG$. Then in fact, when the traces
exist,
$$
\Tr_{\rnp}G=\Tr_{\Bbb R^{n-1}}\tr_n G,\tag 2.18
$$
and there is a similar rule for the operators carried over to the
manifold situation, when the symbol-kernel of $G$ is supported in
$X_c$ and the product measure is used on $X_c$:$$
\Tr_XG=\Tr_{X'}\tr_nG.\tag2.19
$$
In this way, the calculation of traces of s.g.o.s is reduced to the
calculation of traces of $\psi $do's on $X'$, for which we have the results of
[GS95] for operators with symbols in the spaces $S^{m,d,0}(\Gamma )$.
(When $G$ is given as a finite sum of compositions of Poisson and trace
operators, $G=\sum_{j\le J}K_jT_j$, one has by linearity and circular
perturbation 
that $\Tr_XG=\Tr_{X'}(\sum_{j\le J}T_jK_j$), which is a closely related
``reduction to the boundary'' that avoids explicit mention of normal
and tangential variables.) 

Let $\zeta (x_n)$ be a $C^\infty $ function on $\Bbb R$
such that $$
\zeta (x_n) =1\text{ for }|x_n|\le \tfrac13,\;\zeta (x_n) 
\in [0,1]\text{ for } |x_n|\in [\tfrac13,\tfrac23],\;
\zeta (x_n)=0 \text{ for }|x_n|\ge
\tfrac23;\tag2.20
$$
denote $\zeta (x_n/\varepsilon )$ by $\zeta _\varepsilon $,
$\varepsilon >0$. 
Recall from [G01, Lemma 7.1ff.]:

\proclaim{Lemma 2.9}
For a singular Green operator $G$ of class {\rm 0} in the calculus,
$(1-\zeta _\varepsilon )G$ and $G(1-\zeta _\varepsilon )$ have  
symbol-kernels in 
$\Cal S^{-\infty ,-\infty ,-\infty }(\Gamma ,\Cal S_{++})$; hence they are trace-class with
traces that are $O(|\lambda |^{-N})$ for $|\lambda |\to \infty $
in $ \Gamma $, any
$N$.
\endproclaim  

This relies on the fact that such operators
can be written with a factor
$x_n^k$ resp\. $y_n^k$ for any $k$, where Lemma 2.3 applies.

\proclaim{Theorem 2.10} Let $m$, $d$ and $s\in\Bbb Z$, with $s\le
0$.

{\rm (i)} Let $G$ be a $\mu $-dependent singular Green operator of
class $0$, with  polyhomogeneous symbol-kernel in $\Cal
S^{m,d,s-1}(\Gamma ,\Cal S_{++})$ in local coordinates, holomorphic
in $\mu $. If $m+s>
-n$, assume furthermore that the homogeneous terms in the
symbol of $S=\tr_nG$ of degree $m+d+s-j$ with $m+s-j> -n$ are
integrable i $\xi '$. Then $G$ is trace-class and its trace has an asymptotic
expansion in $\mu $ for $|\mu |\to \infty $ in $\Gamma $:
$$
\Tr G\sim \sum_{j\in\Bbb N}c_j\mu ^{m+d+s+n-1-j}+\sum_{k\in\Bbb
N}(c'_k\log\mu +c''_k)\mu ^{d+s-k};\tag2.21
$$
here the coefficients $c_j$ and $c'_k$ with $k=-m+j-n+1$ are
determined from the $j$'th
homogeneous term in the symbol of $\tr_nG$ (are
``local''), whereas the 
$c''_k$ depend on the full operator (are ``global'').
Such an expansion likewise hold for $\Tr_{X'}S$, when $S$ is a $\psi
$do on $X'$ with polyhomogeneous symbol in $S^{m,d,s}(\Gamma )$,
under the same additional assumption as above when $m+s>-n$. 

{\rm (ii)} In particular, if the symbol of $G$ is strongly
polyhomogeneous, then
$$
\Tr G\sim \sum_{j\in\Bbb N}d_j\mu ^{m+d+s+n-1-j}
,\tag2.22
$$
where $d_j$ is determined from the $j$'th homogeneous term in the
symbol of $G$. 
A similar statement holds for $\Tr_{X'}S$ when the symbol of $S$ is strongly polyhomogeneous.
\endproclaim 

\demo{Proof}
This is shown in [G01, Th\. 7.3], but since we need to refer to
specific coefficients, we recall some ingredients of the proof here.
 
By Lemma 2.9, the expansion (2.21) (or (2.22)) is
unaffected by replacing the given s.g.o.\ $G$ by an operator $G'=\zeta
_\varepsilon G\zeta _\varepsilon $ with
symbol-kernel supported in $X_c$, where (2.19) can be used.

We have from (2.12) that the symbol-kernel $\tilde g$ of $G'$ is in$$
\Cal
S^{m,d,s-1}(\Gamma ,\Cal S_{++})\subset \Cal
S^{m+s,d,-1}(\Gamma ,\Cal S_{++})\cap \Cal
S^{m,d+s,-1}(\Gamma ,\Cal S_{++}),
$$
in local coordinates.
Hence, by Proposition 2.8, $\tr_n
\tilde g\in 
S^{m+s,d,0}(\Gamma )\cap S^{m,d+s,0}(\Gamma )$.
Denote $\tr_n G'=S$, then it is a $\psi $do on $X'$ with symbol 
$$s(x',\xi ',\mu )\in 
S^{m+s,d,0}(\Gamma )\cap S^{m,d+s,0}(\Gamma ),\tag2.23$$ 
in local coordinates.
Now one simply applies [GS95, Th\. 2.1] to $S$, and from here on, the
considerations apply to any $\psi $do $S$ on $X'$ satisfying the
stated assumptions:

If $m+s\le -n$,
the inclusion in the first space in (2.23) 
assures trace-class and
integrability in $\xi '$ of all terms in the symbol. An application
of [GS95, Th\. 2.1] then gives after integration in $x'\in X'$ 
that $\Tr_{X'}S$ has an expansion
as in (2.21) but with $d+s$ replaced by $d$ in the second series.
The inclusion in the second space in (2.23) allows us to replace $d$ by the
lower integer $d+s$ in the second
series. To explain the central part of the proof, let $s(x',\xi ',\mu
)\in S^{m',d',0}(\Gamma )$ and let $s'=\mu ^{-d'}s$ (in a localized
situation). Here the $j$'th
homogeneous term $s'_j(x',\xi ',\mu )$ in the symbol $s'(x',\xi ',\mu
)$  has homogeneity degree $m'-j$, denoted $m_j$ in [GS95, Th\.
2.1]. Its contribution to the diagonal value of the kernel $K(x',y',\mu
)$ of $S$ is$$
K_{s'_j}(x',x',\mu )=\int_{\Bbb R^{n-1}}s'_j(x',\xi ',\mu )\,\d\xi '=\Bigl(\int_{|\xi '|\ge |\mu
|}+\int_{|\xi '|\le 1}+\int_{1\le |\xi '|\le |\mu
|}\Bigr)s'_j\,\d\xi ';\tag2.24
$$
here $\int_{|\xi '|\ge |\mu
|}s'_j\,\d\xi '$ gives a power term $c(x')\mu ^{m'-j+n-1}$ (by
homogeneity), $\int_{|\xi 
'|\le 1}s'_j\,\d\xi 
'$ gives a series of power terms, and $\int_{1\le |\xi '|\le |\mu
|}s'_j\,\d\xi '$ gives a log-power term $c'(x')\mu ^{m'-j+n-1}\log \mu
$ if $m'-j+n-1$ is a nonpositive integer and besides this
some power terms. Thus in the original symbol, $s_j=\mu ^{d'}s'_j$ gives a
log-power term   $c'(x')\mu ^{m'+d'-j+n-1}\log \mu
$ when $j\ge m'+n-1$. We note that the log-power terms begin with $c\mu
^{d'}\log \mu $; also the nonlocal terms begin at this power.

If $m+s>-n$, the supplementary
integrability assumption for the terms with $m+s-j>-n$ assures
trace-class, and all terms are treated as above.

In (ii), one gets the refined expansion (2.22) since in the
strongly polyhomogeneous case, the 
homogeneous terms in the symbol are strictly homogeneous
in $(\xi ',\mu )\in \Bbb R^{n-1}\times \{\mu \in\Gamma \mid |\mu |\ge
\varepsilon \}$ and integrable in $\xi '\in\Bbb R^{n-1}$, so that
the terms $d_j\mu ^{m+d+s+n-1-j}$ are produced directly by
integration of the $j$'th symbol in $\xi '$ and $x'$  using the
homogeneity. (Here one does not need to decompose the integral into  
three regions as in (2.24).) 

One can also derive (2.22) from the fact that the strongly polyhomogeneous
symbols are as in [G96] with 
regularity number $\nu =+\infty $, so that full trace expansions with
purely power terms hold as shown there (and recalled in [G01, Prop\.
7.2]).   \qed

\enddemo 

When the operator acts on the sections of a vector bundle over $X$,
one takes the fiber trace in (2.24).
There is a similar result for $\Tr_X P_+$ when $P$ has symbol in
$S^{m,d,s}(\Gamma )$, only with $n-1$ replaced by $n$.

Let us now show how the elementary operators introduced in
(1.18)--(1.20) fit into the calculus.
In the statements below, it is tacitly understood that the symbols
are $N\stimes N$-matrix valued. This could be indicated by adding
$\otimes \Cal L(\Bbb C^N)$ to all the mentioned symbol
spaces, but 
that would make the reading unnecessarily heavy.

\proclaim{Proposition 2.11} 
Consider the operators from Section {\rm 1} as parametrized by $$
\mu =(-\lambda
)^{\frac12}, \text{ in } \Gamma =\{\mu \in\Bbb C\mid \operatorname{Re}\mu
>0\}.\tag2.25
$$

{\rm (i)} The $\psi $do $Q^0_\lambda =(A^2+D_{x_n}^2+\mu ^2)^{-1}$
has symbol in $S^{0,0,-2}_{\operatorname{sphg,ut}}(U\stimes\Bbb R^n,\Gamma )$,
in local trivializations.

{\rm (ii)} The $\psi $do's $$
\aligned
A_\lambda &=(A^2-\lambda )^{\frac12}=(A^2+\mu
^2)^{\frac12}\\
A_\lambda ^{-1}&=(A^2-\lambda )^{-\frac12}=(A^2+\mu
^2)^{-\frac12}
\endaligned\tag2.26
$$ have strongly polyhomogeneous symbols in $S^{0,0,1}(U'\stimes\Bbb
R^{n-1},\Gamma )$ 
resp\. $S^{0,0,-1}(U'\stimes\Bbb R^{n-1},\Gamma )$, in local trivializations. Moreover, for
$m\in\Bbb N$,  
$$\aligned
\partial _ \lambda ^rA_\lambda &\text{ has symbol in }S^{0,0,1-2r}(U'\stimes\Bbb R^{n-1},\Gamma ),\\
\partial _ \lambda ^rA_\lambda^{-1} &\text{ has symbol in
}S^{0,0,-1-2r}(U'\stimes\Bbb R^{n-1},\Gamma ). 
\endaligned\tag2.27$$

{\rm (iii)} The $\psi $do $(A_\lambda +|A|)^{-1}$
has symbol in 
$S^{0,0,-1}(U'\stimes\Bbb R^{n-1},\Gamma )$, in local
trivializations, and its $\partial _ \lambda ^r$-derivatives have symbols in $S^{0,0,-1-2r}(U'\stimes\Bbb R^{n-1},\Gamma )$.

{\rm (iv)} The Poisson operator $K_{A_\lambda }$ and the trace
operator $T_{A_\lambda }$ have strongly
polyhomogeneous symbol-kernels in $\Cal S^{0,0,-1}(U'\stimes\Bbb
R^{n-1},\Gamma ,\Cal S_+) 
$, in local trivializations. Moreover,
$$
\partial _ \lambda ^rK_{A_\lambda }\text{ and }\partial _\lambda
^rT_{A_\lambda }
\text{ have symbol-kernels in 
}S^{0,0,-1-2r}(U'\stimes\Bbb R^{n-1},\Gamma ,\Cal S_+).\tag2.28
$$

\endproclaim 

\demo{Proof} 
(i), essentially known from [GS95], follows from the fact that
$A^2+D_{x_n}^2+e^{2i\theta 
}D_{x_{n+1}}^2$ is an elliptic differential operator on $X'\times\Bbb
R^2$, for any $|\theta |<\frac\pi 2$; then we can for each ray
in $\Gamma $ use Theorem 2.6
(ii) and the fact that elliptic differential operators and their
parametrices satisfy the transmission condition. 

The results of (ii) and (iii) are also essentially known from [GS95]. In (ii),
we can compare $A_\lambda $ and $A_\lambda ^{-1}$ with
$(A^2+e^{2i\theta }D_{x_{n}}^2)^{\frac12}$ resp\. $(A^2+e^{2i\theta
}D_{x_{n}}^2)^{-\frac12}$ defined according to Seeley [S69]; the
latter are elliptic of degree $1$ resp\. $-1$. Then Theorem 2.6 (i)
gives the statement for $r=0$.
The cases $r>0$ are included by use of the formula $\partial _\lambda
A_\lambda =-\frac12 
A_\lambda ^{-1}$ and the composition rules.

Concerning (iii), it is shown in the proof of [GS95, Prop\. 3.5] that
$A_\lambda (|A|+A_\lambda )^{-1}$ has symbol in $S^{0,0}(\Gamma )$,
equal to $S^{0,0,0}(\Gamma )$ by (2.11), so
the result follows for $r=0$ by composition with $A_\lambda ^{-1}$,
using (ii) and the composition rule (xiii) in Theorem 2.7. For $r>0$,
we use that $$ 
\partial _ \lambda ^r(|A|+A_\lambda )^{-1}=\sum_{j+k=2m, j,k\ge 1}c_{jk}(|A|+A_\lambda
)^{-1-j}A_\lambda ^{-k},\tag2.29
$$
so that (iii) follows from the previous results by use of 
rule (xiii).

For (iv), we apply Theorem 2.6 (iii). In fact, when
$\mu =e^{i\theta }\varrho $ is replaced by $e^{i\theta }\xi 
_{n+1}$, the symbol-kernel of $K_{A_\lambda }$ is replaced by the
symbol-kernel  of
the solution operator (in a parametrix sense) $\widetilde K\:\varphi
\mapsto u$ of the
Dirichlet problem$$ 
\aligned
(A+D_{x_n}^2+e^{i2\theta }D^2_{x_{n+1}})u(x',x_n,x_{n+1})&=0\text{ on
}\Bbb R^{n-1}\stimes \Bbb R_+\stimes\Bbb R,\\
u(x',0, x_{n+1})&=\varphi (x',x_{n+1}) \text{ on }\Bbb
R^{n-1}\stimes\Bbb R,
\endaligned\tag2.30
$$
which is a standard Poisson operator relative to $\Bbb
R^{n-1}\stimes\Bbb R_+\stimes\Bbb R$.
Then Theorem 2.6 (iii) implies that the symbol-kernel of $K_{A_\lambda
}$ is in $\Cal S^{0,0,-1}(\Gamma ,\Cal S_+)$. The $\lambda
$-derivatives are included by functional calculus (cf\.
(1.19)) and composition rules. There is a similar proof
for $T_{A_\lambda }$. 
\qed

\enddemo

\subhead 3. Preservation of log-terms under general perturbations 
\endsubhead

We shall here study the resolvent $R_\lambda =(\Delta _B-\lambda
)^{-1}$
by use of the representation (1.14).
It is shown in [GS95, Th\. 3.9] for the non-product case (with
$\Pi =\Pi _>+B_0$), in
[G99, Cor\. 8.3] for the general case, that $\Cal R_\mu $ has the structure$$
\Cal R_\mu =\Cal Q_{\mu ,+}+\Cal G_\mu ,\quad \Cal G_\mu =\Cal K_\mu \Cal S_\mu \Cal T_\mu ,\tag3.1
$$
where $\Cal S_\mu $ is a weakly polyhomogeneous 
$\psi $do on $X'$ with symbol in $S^{0,0}(\Gamma )$, $\Cal K_\mu $ is a
strongly polyhomogeneous Poisson operator of degree $-1$, and
$\Cal T_\mu $ is a strongly polyhomogeneous trace operator of class 0 and
degree $-1$.  $\Cal Q_\mu $ is the parametrix described in (1.12).

From the point of view of
the more recent calculus in recalled in Section 2, $\Cal
S_\mu $ has symbol 
in $S^{0,0,0}(\Gamma )$ (cf\. (2.11)), and $\Cal K_\mu $ 
and $\Cal T_\mu $ have symbol-kernels in $\Cal S^{0,0,-1}(\Gamma
,\Cal S_+)$ since they are strongly polyhomogeneous (cf\. Theorem 2.6).
So by the elementary composition rules (xii) and (vii)
in Theorem 2.7, we find that $\Cal G_\mu $ has symbol-kernel in 
$\Cal S^{0,0,-2}(\Gamma ,\Cal
S_{++})$, in local trivializations.

We use this to see from (1.14), (1.15) that 
$$
R_{-\mu ^2}=Q_{-\mu ^2,+}+G_{-\mu ^2},\quad G_{-\mu^2} =\mu
^{-1}\pmatrix 1&0   
\endpmatrix \Cal G _\mu \pmatrix 1\\0
\endpmatrix   ,\tag3.2
$$
where 
$G_{-\mu ^2} $ has symbol-kernel in $\Cal S^{0,-1,-2}(\Gamma ,\Cal
S_{++})$, in local trivializations.
$Q_{-\mu ^2}$ is strong\-ly polyhomogeneous of degree $-2$ and has symbol in
$S^{0,0,-2}_{\operatorname{spgh,ut}}(\Gamma )$, by Theorem 2.6 (ii).

\proclaim{Lemma 3.1}
For any $r\in\Bbb N$, 
$$
 R_{-\mu ^2}^r=(Q^r_{-\mu ^2})_+ +G_{-\mu ^2} ^{(r)},\tag3.3
$$
where $Q_{-\mu ^2}^r$ is strongly polyhomogeneous of degree $-2r$,
with symbol in $S^{0,0,-2r}_{\operatorname{spgh,ut}}(\Gamma )$,
 and $G_{-\mu ^2} ^{(r)}$ has symbol in $\Cal S^{0,-r,-1-r}(\Gamma ,\Cal
S_{++})$, in local trivializations.
\endproclaim 

\demo{Proof}
The statement for
the case $r=1$ is shown above. The iterated expressions:$$
R^r_{-\mu ^2}=(Q_{-\mu ^2,+}+G_{-\mu ^2} )\circ\dots\circ(Q_{-\mu ^2,+}+G_{-\mu ^2} ).
\tag3.4$$
are included by use of rules (iv)--(vi) and (xiv) in Theorem 2.7.
\qed\enddemo 

\example{Remark 3.2}
By a direct study of the resolvent of $D_\Pi ^*D_\Pi $, as carried
out for the non-product case with $\Pi =\Pi _\ge$ in [G92], and
for more general cases in [G02], one can show that $G_{-\mu ^2} =K_\mu S_\mu
T_\mu $, where $K_\mu $ and 
$T_\mu $ are strongly polyhomogeneous of degree $-1$ and $S_\mu $ has
symbol is in $S^{0,0,-1}(\Gamma )$; hence $G_{-\mu ^2} $ has
symbol-kernel in $\Cal S^{0,0,-3}(\Gamma ,\Cal S_{++})$, which leads
to the conclusion that $G_{-\mu ^2} ^{(r)}$ in fact
has symbol-kernel in $\Cal S^{0,0,-1-2r}(\Gamma ,\Cal S_{++})$. However, the
above information suffices for the results on perturbations that we
pursue here.
\endexample

It is well-known how $\Tr Q^r_{-\mu ^2,+}$ has an asymptotic
development in pure powers of $\mu $ (when $2r>n$). When we in addition
apply Theorem 2.10 to $G^{(r)}_{-\mu ^2}$ we get (1.8), reconfirming the
result of [GS95]. 
Let us also describe the trace expansion in cases where $R^r_{-\mu
^2}$ is composed with a differential operator:

\proclaim{Theorem 3.3} Let $F$ be a differential operator of order
$m'$. Then for $r>\frac{n+m'}2$, $FR^r_{-\mu ^2}$ is trace-class and
has an expansion$$
\Tr F R^r_{-\mu ^2 }
\sim
\sum_{-n\le k<0} \tilde a_{ k}\mu ^{{{m'} -k}-2r}+ 
\sum_{k\ge 0}\bigl({ \tilde a'_{ k}}\log \mu +{\tilde a''_{
k}}\bigr)\mu ^{{{m'} -k}-2r}.\tag3.5
$$
If $F$ is tangential on $X_c$, then
$$
\Tr F R^r_{-\mu ^2 }
\sim
\sum_{-n\le k<m'} \tilde a_{ k}\mu ^{{{m'} -k}-2r}+ 
\sum_{k\ge m'}\bigl({ \tilde a'_{ k}}\log \mu +{\tilde a''_{
k}}\bigr)\mu ^{{{m'} -k}-2r}.\tag3.6
$$
The coefficients $\tilde a_k$ and $\tilde a'_k$ are locally determined. If $m'$ is odd,
$\tilde a_{-n}=0$.
\endproclaim 

\demo{Proof} Since the operator is of order $m'-2r$ for fixed $\mu $,
it is trace-class when \linebreak $m'-2r<-n$, i.e., $r>\frac{n+m'}2$. It is
well-known that the $\psi
$do part $FQ^r_{-\mu ^2,+}$ has an expansion $\sum_{k\ge -n}c_k\mu
^{m'-2r-k}$
without logarithmic terms. Here the terms with $k$ even/odd vanish when
$m'$ is odd/even, respectively, since they are defined by
integration in $\xi \in \Bbb R^n$ of symbols that are odd in $\xi $.
One finds using Lemma 2.3 that the s.g.o\. 
part $FG^{(r)}_{-\mu ^2}$ has 
symbol-kernel in $\Cal S^{0,-r,m'-r-1}(\Gamma ,\Cal S_{++})$ in local
trivializations, so it follows from Theorem 2.10, applied with $
(m,d,s)=(0,-r,m'-r)$,
that this term contributes a trace expansion as in
(3.5) but starting with $k=1-n$. Then (3.5) follows by summation;
in particular, the term with $k=-n$ vanishes if $m'$ is odd.

If $F$ is tangential, the symbol-kernel of $FG^{(r)}_{-\mu ^2}$ is
instead in $\Cal 
S^{m',-r,-r-1}(\Gamma ,\Cal S_{++})$, so we can use Theorem
2.10 with $(m,d,s)=(m',-r,-r)$, which lowers the starting power $d+s$ in the
series with logarithms to $-2r$; this results in (3.6).\qed
\enddemo

Now we consider two choices of $D$ as in (1.1), with the same
well-posed choice
of boundary condition:$$
D_i=\sigma (\partial _{x_n}+A_{1i}),\quad \Pi \gamma _0u=0,
\quad i=1,2.
\tag3.7
$$
The associated other operators in the direct and the doubled-up
situation will be marked by index $1$ or $2$. The realizations of
$\Cal D_1$ and $\Cal D_2$ are defined by the same boundary condition
(1.11), so the full systems in the doubled-up situations are:$$
\pmatrix \Cal D_1+\mu \\ \Cal B\gamma _0 
\endpmatrix ,\quad
\pmatrix \Cal D_2+\mu \\ \Cal B\gamma _0 
\endpmatrix ,\tag3.8
$$
with inverses$$
\pmatrix \Cal D_1+\mu \\ \Cal B\gamma _0 
\endpmatrix ^{-1}=\pmatrix \Cal R_{1,\mu }&\Cal K_{1,\mu } 
\endpmatrix ,\quad
\pmatrix \Cal D_2+\mu \\ \Cal B\gamma _0 
\endpmatrix = \pmatrix \Cal R_{2,\mu }&\Cal K_{2,\mu } 
\endpmatrix. \tag3.9
$$
We can write (on $X_c$)
$$D_1-D_2=x_n^l\overline P_l,\quad D^*_2-D^*_1=x_n^l\overline P^*_l\tag3.10
$$
for some $l\ge 0$, some tangential $x_n$-dependent first-order
differential operator $\overline P_l$.
Then
$$
 \Cal D_2-\Cal D_1=x_n^l\overline {\Cal
P}_l,\text{ where }\overline {\Cal P}_l=\pmatrix  0&-\overline P_l^*\\ \overline P_l&0
\endpmatrix .\tag3.11  
$$
Let us (somewhat abusively) apply the notation $x_n^l\overline {\Cal
P}_l$ to $\Cal D_2-\Cal D_1$ on all of $X$. Then
$$
\aligned
\left(\matrix \Cal R_{2, \mu}& \Cal K_{2, \mu}\endmatrix\right)
&=\left(\matrix \Cal R_{1, \mu}& \Cal K_{1, \mu}\endmatrix\right)
\left(\matrix \Cal D_1+\mu \\
\Cal B\gamma _0 \endmatrix\right)
\left(\matrix \Cal R_{2, \mu}& \Cal K_{2, \mu}\endmatrix\right)\\
&=\left(\matrix \Cal R_{1, \mu}& \Cal K_{1, \mu}\endmatrix\right)
\left(\matrix \Cal D_2-x_n^l\overline{\Cal P}_l+\mu \\
\Cal B\gamma _0 \endmatrix\right)
\left(\matrix \Cal R_{2, \mu}& \Cal K_{2, \mu}\endmatrix\right)\\
&=\left(\matrix \Cal R_{1, \mu}& \Cal K_{1, \mu}\endmatrix\right)
-\left(\matrix \Cal R_{1,\mu }x_n^l\overline{\Cal P}_l\Cal R_{2,\mu }&\Cal
R_{1,\mu }x_n^l\overline{\Cal P}_l\Cal
K_{2,\mu }\endmatrix\right).\endaligned\tag3.12
$$
In particular,$$
\Cal R_{2,\mu }-\Cal R_{1,\mu }=-\Cal R_{1,\mu }x_n^l\overline{\Cal P}_l\Cal
R_{2,\mu }.\tag3.13 
$$

\proclaim{Theorem 3.4} Consider $D_{1,\Pi }$ and $D_{2,\Pi }$ defined
from {\rm (3.7)}. 
For any $r\in\Bbb N$, write$$
R^r_{i,-\mu ^2}=(Q^r_{i,-\mu ^2})_++G^{(r)}_{i,-\mu^2},\;i=1,2;\tag3.14
$$
according to Lemma {\rm 3.1}, $Q^r_{i,-\mu ^2}$ has symbol in
$S^{0,0,-2r}_{\operatorname{spgh,ut}}(\Gamma )$ and $G^{(r)}_{i,-\mu ^2}$
has symbol-kernel in $\Cal S^{0,-r,-1-r}(\Gamma ,\Cal S_{++})$ in local
trivializations. Write
$$
R_{2,-\mu ^2}^r-R_{1,-\mu ^2}^r=(Q^r_{2,-\mu ^2}-Q^r_{1,-\mu ^2})_++
\overline G^{(r)}_\mu .\tag3.15
$$
When {\rm (3.10)} holds for some some
$l\ge 0$
(with a first-order tangential differential operator $\overline P_l$), then
$\overline G^{(r)}_{\mu
}$
has symbol-kernel in $\Cal S^{1,-r,-2-r-l}(\Gamma ,\Cal S_{++})$, in local
trivializations.
\endproclaim 

\demo{Proof} We begin with the case $r=1$. 
By (3.13) and (3.2),
$$\multline
\Cal R_{2,\mu }-\Cal R_{1,\mu }=-\Cal R_{1,\mu }x_n^l\overline{\Cal P}_l\Cal
R_{2,\mu }
=-(\Cal Q_{1,\mu,+}+\Cal G_{1,\mu })x_n^l\overline{\Cal P}_l (\Cal
Q_{2,\mu,+}+\Cal G_{2,\mu })\\
=-\Cal Q_{1,\mu,+}x_n^l\overline{\Cal P}_l \Cal Q_{2,\mu,+}-\Cal
Q_{1,\mu,+}x_n^l\overline{\Cal P}_l \Cal G_{2,\mu }-
\Cal G_{1,\mu }x_n^l\overline{\Cal P}_l \Cal Q_{2,\mu,+}-
\Cal G_{1,\mu }x_n^l\overline{\Cal P}_l \Cal G_{2,\mu }
 .\endmultline\tag3.16 
$$
We first show how the three last terms
are treated by Theorem 2.7 and Lemma 2.3. The lemma shows that
multiplication of a singular 
Green symbol-kernel by $x_n^l$ or $y_n^l$ lowers the third upper index by $l$ steps.
Thus $x_n^l\Cal G_{2,\mu }$ 
has symbol-kernel in $\Cal S^{0,0,-2-l}(\Gamma ,\Cal S_{++})$, in local
trivializations. A similar result holds for $\Cal G_{1,\mu }x_n^l$, where
the composition with $x_n^l$ to the right has the effect of
multiplying the symbol-kernel with $y_n^l$. 
The composition with the  $\mu
$-independent first-order tangential  differential operator $\Cal
P_l$ has the effect of lifting the first upper index by 1 step.
When we take these
effects into account and use the composition rules, we find that the
three last terms in (3.16) are s.g.o.s with symbol-kernels in 
$\Cal S^{1,0,-3-l}(\Gamma ,\Cal
S_{++})$, in local trivializations.
The remaining term equals
$$
\Cal Q_{1,\mu,+}x_n^l\overline{\Cal P}_l \Cal Q_{2,\mu,+}
=(\Cal Q_{1,\mu }x_n^l\overline{\Cal P}_l\Cal Q_{2,\mu })_+-G^+(\Cal Q_{1,\mu
})x_n^l\overline{\Cal P}_l G^-(\Cal Q_{2,\mu }),\tag3.17
$$
cf\. (xiv) of Theorem 2.7. Here the $\psi $do part is strongly
polyhomogeneous of order $-1$ 
and the $G^{\pm}(\Cal Q_{i,\mu })$ 
have symbol-kernels in
 $\Cal S^{0,0,-2}(\Gamma ,\Cal S_{++})$, in local trivializations. By Lemma
2.3, $G^{+}(\Cal Q_{i,\mu })x_n^l$ 
has symbol-kernel in
 $\Cal S^{0,0,-2-l}(\Gamma ,\Cal S_{++})$. $\overline{\Cal P}_l$
lifts the first upper index by 1. Then by
Theorem 2.7 (vi), the last term in (3.17) has
symbol-kernel in  $\Cal S^{1,0,-3-l}(\Gamma ,\Cal
S_{++})$, in local trivializations. 

Now by (1.14),
$$
R_{2,-\mu ^2}-R_{1,-\mu ^2}=\mu ^{-1}\pmatrix 1&0 
\endpmatrix (\Cal R_{2,\mu }-\Cal R_{1,\mu }) \pmatrix 1\\0
\endpmatrix  
=-\mu ^{-1}\pmatrix 1&0 
\endpmatrix \Cal R_{1,\mu }x_n^l\overline{\Cal P}_l\Cal
R_{2,\mu } \pmatrix 1\\0
\endpmatrix,
\tag3.18 
$$
so from what we showed for (3.16) follows in view
of (3.2) that 
$R_{2,-\mu ^2}-R_{1,-\mu ^2}$ is the sum of a truncated
$\psi $do with symbol in $S^{0,0,-2}_{\operatorname{sphg,ut}}(\Gamma
)$ and an s.g.o\. with symbol-kernel in $\Cal S^{1,-1,-3-l}(\Gamma ,\Cal
S_{++})$, in local trivializations. This shows the assertion for $r=1$.

For higher $r$, we use the formula$$
R_{2,-\mu ^2}^r-R_{1,-\mu ^2}^r=
(R_{2,-\mu ^2}-R_{1,-\mu ^2})(R_{2,-\mu ^2}^{r-1}+R_{2,-\mu
^2}^{r-2}R_{1,-\mu ^2}+\cdots +R_{1,-\mu ^2}^{r-1}).\tag3.19
$$
The first factor is described above. For
the terms in the second factor we use the information in (3.14)ff. An
application of the 
composition rules gives an operator whose $\psi $do part has symbol
in $S^{0,0,-2r}_{\operatorname{sphg,ut}}(\Gamma )$ and whose s.g.o\.
part has symbol-kernel in$$
\Cal S^{1,-1,-3-l}(\Gamma ,\Cal S_{++})\circ
\Cal S^{0,-(r-1),-1-(r-1)}(\Gamma ,\Cal S_{++})\subset
\Cal S^{1,-r,-2-r-l}(\Gamma ,\Cal S_{++}).\quad\square
$$ 
\enddemo 

It is remarkable in this result (and important for the applications below)
that both factors $x_n^l$ (for $l>0) $ and $\Cal R_{2,\mu }$ 
lower the third upper index, whereas $\overline{\Cal P}_l$ only lifts
the first upper index. In the case $l=0$, the proof in fact allows
$D_2-D_1$ to be an
arbitrary first-order tangential differential operator (as long as
ellipticity is respected).

Pursuing the indications in Remark 3.2, one could get a still better
result, placing the decrease by $2r$ fully in the third upper index. However,
the above result suffices to conclude:

\proclaim{Theorem 3.5} Let $D_1$ and $D_2$ be two first-order
elliptic operators on
$X$ as in {\rm (1.1), (3.7)}, provided with the same well-posed
boundary condition $\Pi \gamma 
_0u=0$ (with $\Pi $ being an orthogonal pseudodifferential projection); then
{\rm (3.10)} holds for some $l\ge 0$, and we denote the largest such
integer by $l$. Let
$D_{1,\Pi }$ and $D_{2,\Pi }$ be the realizations defined by the
boundary condition $\Pi \gamma _0u=0$, and let $\Delta
_{i,B}=D_{i,\Pi }^*D_{i,\Pi }$. Let $F$ be a differential operator in
$E_1$ of 
order ${m'} $ and let $r>\frac {n+{m'} }2$. Then$$
F(R_{2,\lambda }^r-R_{1,\lambda }^r)=(F(Q^r_{2,\lambda
}-Q^r_{1,\lambda }))_++
F\overline G^{(r)}_\mu ,\tag3.20
$$
where $F(Q^r_{2,-\mu ^2}-Q^r_{1,-\mu ^2})$ has symbol in
$S^{0,0,m'-2r}_{\operatorname{spgh,ut}}(\Gamma )$ and $F\overline
G^{(r)}_\mu $has symbol-kernel in $\Cal S^{1,-r,m'-2-r-l}(\Gamma ,\Cal
S_{++})$
(in $\Cal S^{m'+1,-r,-2-r-l}(\Gamma ,\Cal S_{++})$ if $F$ is tangential),
in local trivializations.

The $\psi $do part has an asymptotic trace expansion
$$\Tr [(F (Q^r_{2,\lambda }-Q^r_{1,\lambda }))_+]
\sim
\sum_{-n\le k<\infty  } \tilde p_{ k}(-\lambda )^{\frac{{m'} -k}2-r},
\tag3.21
$$
where $\tilde p_k=0$ for $k-m'+n$ odd.

The s.g.o.\ part has an asymptotic trace expansion
$$
\Tr [F\overline G^{(r)}_\mu)]
\sim
\sum_{-n+1+l\le k< k_0 } \tilde g_{ k}(-\lambda )^{\frac{{m'}
-k}2-r}+ 
\sum_{k\ge k_0}\bigl({ \tilde g'_{ k}}\log (-\lambda
)+{\tilde g''_{ k}}\bigr)(-\lambda )^{\frac 
{{m'} -k}2-r},\tag3.22
$$
where $$
k_0=l+1 \text{ when $F$ is general},
k_0=m'+l+1 \text{ when $F$ is tangential on $X_c$.}
\tag3.23 
$$
It follows that
$$
\Tr [F  (R^r_{2,\lambda }-R^r_{1,\lambda })]
\sim
\sum_{-n\le k< k_0 } \tilde c_{ k}(-\lambda )^{\frac{{m'}
-k}2-r}+ 
\sum_{k\ge k_0}\bigl({ \tilde c'_{ k}}\log (-\lambda
)+{\tilde c''_{ k}}\bigr)(-\lambda )^{\frac 
{{m'} -k}2-r},\tag3.24
$$
with $k_0$ as above. For $k\le l-n$, the $\tilde c_k$ vanish when
$k-m'+n$ is odd. 

The coefficients $\tilde c_k$ and $\tilde c'_k$ are locally determined.

\endproclaim 

\demo{Proof} Recall that $\lambda =-\mu ^2$. 
The statement on the decomposition (3.20) follows from Theorem 3.4
and Lemma 3.1.
The trace expansion of the $\psi
$do part is well-known. The trace
expansion of the s.g.o.\ part is obtained in general by application of
Theorem 2.10, with
$$
m=1,\quad  d=-r , \quad s={m'} -1-r-l.$$
 Here
$m+s={m'} -r-l$, so if  $r +l\ge n+{m'} $, there is no need to check
integrability conditions. We only assume $r>\frac {n+{m'} }2$; this is
allowed because the singular Green part of the resolvent power is in
fact of order $-2r$ (degree 
$-2r-1$) for each fixed 
$\lambda $ (since it comes from the resolvent of an elliptic problem
of order $2r$), hence $tr_n$ of its composition with $F$ is a $\psi $do on $X'$ of order
$m'-2r$. So all the homogeneous terms in the considered $\psi $do
symbol are integrable 
in $\xi '$, when $r>\frac{n+{m'} }2$.
Then Theorem 2.10 gives an expansion of the trace of the s.g.o\. part
of the form (2.21), with $m+d+s={m'} -2r-l$, $d+s={m'} -1-2r-l$
(resp\. $-1-2r-l$ if $F$ is tangential). Here the expansion starts
with the power $m'-l+n-1-2r$, and the log-terms start with the power
$m'-1-l-2r$ (resp.\ $-1-l-2r$ if $F$ is tangential), so (3.22)
is obtained after some relabelling.

When we add the contributions, we find (3.24).
\qed
\enddemo 

\example{Remark 3.6}
Note in particular that the terms with ``global'' coefficients
$\tilde c''_k$ begin with the
power $(-\lambda 
)^{\frac {{m'} -l-1}2-r}$ for general $F$, $(-\lambda 
)^{\frac {-l-1}2-r}$ when $F$ is tangential.
\endexample

There is a similar result for $D_{i,\Pi }{D_{i,\Pi }}^*$. We have furthermore:

\proclaim{Corollary 3.7} Hypotheses and definitions as in Theorem
{\rm 3.5}. There are 
expansions$$\aligned
\Tr [F (e^{-t\Delta _{2,B}}-e^{-t\Delta _{1,B}})]
&\sim
\sum_{-n\le k<k_0 } c_{ k}t^{\frac{k-{m'}}2}+ 
\sum_{k\ge k_0}\bigl({  c'_{ k}}\log t+{ c''_{ k}}\bigr)t^{\frac 
{k-{m'}}2},\\
\Gamma (s)\Tr [F (\Delta _{2,B}^{-s}-\Delta _{1,B}^{-s})]
&\sim
\sum_{-n\le k<k_0 } \frac{c_{ k}}{s+\frac{k-{m'}}2}+
\sum_{k\ge k_0}\Bigl(\frac{-  c'_{ k}}{(s+{\frac 
{k-{m'}}2})^2}+\frac{ c''_{ k}}{s+\frac 
{k-{m'}}2}\Bigr)\\
&\qquad-\frac{\Tr[F(\Pi
_0(\Delta _{2,B})-\Pi _0(\Delta _{1,B}))]}s.
\endaligned\tag3.25$$

The coefficients $ c_k$ and $c'_k$ are locally determined; the $c_k$,
$c'_k$, $c''_k$ are 
proportional to $\tilde c_k$, $\tilde c'_k$, $\tilde c''_k$ in {\rm
(3.24)} by universal factors. For
$k\le l-n$, the $c_k$
vanish when $k-m'+n$ is odd.
\endproclaim 

\demo{Proof} Here one uses the transition formulas explained e.g.\ in
[GS96]. The passage from (3.24) to the zeta function expansion in the
second formula of 
(3.25) is based on Corollary 2.10 there, and the passage to the heat
trace expansion is based on Section 5 there.
The subtracted term in the second line is explained by the fact that
$\Delta _{i,B}^{-s}$ is 
defined to be zero on 
$V_0(\Delta _{i,B})$.
\qed
\enddemo 

We can also formulate the result as follows:

\proclaim{Corollary 3.8} Hypotheses and definitions as in Theorem
{\rm 3.5}. For the 
trace expansion$$
\Tr F R^r_{1,\lambda  }
\sim
\sum_{-n\le k<0} \tilde a_{ k}(-\lambda )^{\frac{{m'} -k}2-2}+ 
\sum_{k\ge 0}\bigl({ \tilde a'_{ k}}\log (-\lambda ) +{\tilde a''_{
k}}\bigr)(-\lambda ) ^{\frac{{m'} -k}2-r}\tag3.26
$$ (the summation limit $0$ replaced by $m'$ if $F$ is tangential), the
 replacement of $D_1$ by
$D_2$ leaves the log-coefficients 
$\tilde a'_k$ invariant for $ k<k_0$. The other coefficients with
$k<k_0$ are modified only by local terms; those with
$ k\le l-n$ and $k-m'+n$ odd are invariant. 

There are similar results for the associated heat trace and zeta function.
\endproclaim 
\demo{Proof} (3.26) is a reformulation of (3.5)--(3.6). Since 
$F R^r_{2,\lambda  }=F R^r_{1,\lambda  }+F (R^r_{2,\lambda
}-R^r_{1,\lambda  })$, the result follows for $F R^r_{2,\lambda  }$
by addition of (3.24) to (3.26).\qed
\enddemo 

Recall from [GS96]
that when $F=\varphi ^0$ (a morphism 
independent of $x_n$ on $X_c$), the coefficient
$a'_1$ in (0.1) and (1.9) in the product case, with $\Pi $ equal to $\Pi _>$
plus a projection in the nullspace $V_0(A)$, satisfies:
$$
a'_1=-\pi ^{-1}e_1(\varphi ^0,A^2),\tag 3.27
$$ where $e_1$
is the coefficient of $t^{\frac12}$ in the heat trace expansion for $A^2$
on $X'$:
$$
{\Tr}(\varphi ^0e^{-tA^2})\sim \sum_{k= 1-n}^\infty  e_k(\varphi
^0,A^2)t^{\frac k2}.
$$
Here $e_k(\varphi ^0,A^2)=0$ for $k-n+1$ odd; in particular, $e_1=0$
if $n$ is odd. By [G01$'$], the value $a'_1$ is the same also for
projections $\Pi =\Pi _>+\Cal S$ with $\Cal S$ of order $\le -n-1$.

The above methods moreover allow us to conclude:

\proclaim{Theorem 3.9} 
The coefficient $a'_1$ in {\rm (0.1)} and 
{\rm (1.9)} (as well as the coefficient $\tilde a'_1$ in {\rm (1.8)}), is 
the same for a non-product type operator $D$
 {\rm (1.1), (1.2)} over $X_c$ with volume form $v(x)\,dx$ and the
associated
product type operator $D^0$
{\rm (1.3)} with volume form $v(x',0)\,dx$, when $P_0=0$ and
$\partial _{x_n}v(x',0)=0$.
Moreover,
the coefficient $a''_1$ (as well as $\tilde a''_1$) differs in the
cases of $D$ and $D^0$
by a local contribution only.  
Here, when $\varphi =\varphi ^0$ (independent of $x_n$ on $X_c$) and $\Pi
=\Pi _>+\Cal S$, $\Cal S$ of order $\le -n-1$, $a'_1$
satisfies {\rm (3.27)}. 

Generally,  when $D=D^0+x_n^l\overline P_l$ on $X_c$ for some $l$,
and $\partial _{x_n}^jv(x',0)=0$ for $1\le j\le l$, then the terms
$a'_k$ for $0\le k\le l$ are the same in the expansions for $D$ and
for $D^0$, and the nonlocal terms $a''_k$ differ by local
contributions only.
\endproclaim

\demo{Proof} Consider the first mentioned case, where $l=1$; here
$D=D^0+x_n\overline P_1$ with $\overline P_1=\sigma P_1$. Note
that the adjoints $D^*$ of $D$ and ${D^0}'$ of $D^0$ are defined
differently because of 
the different volume forms. However, as noted after (1.2), the
hypothesis $\partial _{x_n}v(x',0)$ assures that
$D^*$ has the form $(-\partial _{x_n}+A+x_nP'_1)\sigma
^*={D^0}'+x_n\overline P'_1$ with  $\overline P'_1=P'_{1}\sigma ^*$.
Then on $X_c$,$$
\Cal D -\Cal D^0= x_n\overline {\Cal
P}_1,\text{ where }\overline {\Cal P}_1=\pmatrix  0&-\overline P'_1\\ \overline P_1&0
\endpmatrix .
$$
The proofs of Theorems 3.4 and 3.5 extend immediately to this situation.
In the expansion corresponding to (3.22) in this case, the second sum
begins at the index $k_0=2$.
This shows the stability, and the statement on the value
follows from the remarks before the theorem.

In the case with general $l$, the hypotheses assure (in view of the
remarks after (1.2)) that $\Cal D -\Cal D^0=x_n^l\overline{\Cal P}_l$
for a suitable first-order tangential operator $\overline{\Cal P}_l$,
and the proof goes as in Theorems 3.4 and 3.5.
\qed
\enddemo

We recall (e.g\. from [GG98]) that $e_1$ is 
generically nonzero when $n$ is even (cf\. also (3.29) below).
 
Note also the general result in the case $l=0$, where $P_0$ can be
nonzero: Here $a'_0$ is preserved when $D$ is replaced by $D^0$, and
$a''_0$ is perturbed only by local 
contributions.
This stability was established for the case with
$\Pi =\Pi _\ge(A)$ in [G92], and in
[GS95] for $\Pi =\Pi _>(A)$ plus certain finite rank projections.

But a perturbation with $P_0\ne 0$
will in general change $a'_1$:

\example {Remark 3.10} Consider the simple case where $P_0=\alpha I$
(on $X_c$),
$\alpha \in\Bbb R$, and $F=\varphi ^0$ as above. 
If we also change the boundary projection to $\Pi _\ge (A+\alpha I)$, we
have a new product case, where the operator is $$
D^0_\alpha =\sigma
(\partial _{x_n}+A+\alpha I) \text{ (on $X_c$) with boundary condition }
\Pi _{\ge}(A+\alpha I)\gamma _0u=0,\tag 3.28
$$ 
and the
heat trace has an expansion as in (0.1) with
log-coefficients that vanish according to the description after (0.1). 
Here the coefficient $
a'_1
$ equals $-\pi ^{-1}e_1(\varphi ^0,(A+\alpha I)^2)$,  by
the preceding explanation.
Now since $\Pi _\ge(A+\alpha )-\Pi (A)$ is a finite linear combination of
eigenprojections of $A$, it is a $\psi $do of order $-\infty $, so a
replacement of $\Pi _\ge(A)$ by \linebreak$\Pi
_\ge(A+\alpha )$ in the boundary condition (1.6) leaves all log-terms
invariant by the results in [G99] (elaborated in [G01$'$]). 
Thus in fact, when we return to the boundary condition (1.6), the term
$a'_1$ in the trace expansion for 
$D^0_\alpha $ is also equal to $-\pi ^{-1}e_1(\varphi ^0,(A+\alpha
)^2)$.

It vanishes for $n$ odd, but let us consider the case $n$ even. As
recalled in [GG98], 
differentiation and comparison of the expansions of
${\Tr}(\varphi ^0e^{-t(A+\alpha )^2})$  and 
${\Tr}(\varphi ^0(A+\alpha )e^{-t(A+\alpha )^2})$ leads to the following formula (with a nonzero
integer factor $m(n) $) $$
\partial _\alpha ^{n}e_1(\varphi ^0,(A+\alpha )^2)=m(n)\,
e_{1-n}(\varphi ^0,(A+\alpha )^2)\ne 0,\tag3.29
$$  
which shows that $a'_1$ {\it for $D^0_\alpha $ with boundary condition} (1.6)
{\it is not constant in $\alpha $ when $n$ is even}. 
\endexample

One can similarly study the differences connected with eta functions,
$$\gathered
\Tr [F \psi  (D _{2,\Pi}R^r_{2,\lambda }-D _{1,\Pi}R^r_{1,\lambda })],\quad
\Tr [F \psi (D _{2,\Pi}e^{-t\Delta _{2,B}}-D _{1,\Pi}e^{-t\Delta 
_{1,B}})]\\
\text{ and } \Gamma (s)\Tr [F \psi (D _{2,\Pi}\Delta 
_{2,B}^{-s}-D _{1,\Pi}\Delta  _{1,B}^{-s})],\endaligned
$$
by considerations as above, departing from the formula inferred from (1.13):
$$
D _{2,\Pi}R_{2,-\mu ^2 }-D _{1,\Pi}R_{1,-\mu ^2}=-\pmatrix 0&1
\endpmatrix (\Cal R_{2,\mu }-\Cal R_{1,\mu })\pmatrix 1\\0
\endpmatrix \tag3.30
$$
(note that $R_{i,\lambda }$ maps into the domain of $D_{i,\Pi})$),
and using again the considerations on the terms in (3.16). Here the
expressions with higher powers are included by use of the formula
$$
\multline
 D_{2}R_{2,\lambda }^r- D_{1}R_{1,\lambda }^r
= (D_{2}R_{2,\lambda  }-D_{1}R_{1,\lambda  })R_{2,\lambda 
}^{r-1}\\ + D_{2}R_{2,\lambda }^r(R_{2,\lambda  }-R_{1,\lambda
})(R_{2,\lambda }^{r-2}+\cdots +R_{1,\lambda }^{r-2})
\endmultline\tag3.31
$$ 
(also used in [G01$'$]). This leads to:

\proclaim{Theorem 3.11}
Hypotheses of Theorem {\rm 3.5}.

{\rm (i)} 
We have for any
$r\ge 1$:
$$
D_{2,\Pi}R_{2,-\mu ^2}^r-D_{1,\Pi}R_{1,-\mu ^2}^r=(D_2Q^r_{2,-\mu
^2}-D_1 Q^r_{1,-\mu ^2})_++
\widetilde G^{(r)}_\mu ,\tag3.32
$$
where $D_2Q^r_{2,-\mu ^2}-D_1Q^r_{1,-\mu ^2}$ has symbol in 
$S^{1,0,-2r}_{\operatorname{spgh,ut}}(\Gamma )$ and $\widetilde G^{(r)}_{\mu
}$
has symbol-kernel in $\Cal S^{1,1-r,-2-r-l}(\Gamma ,\Cal S_{++})
+\Cal S^{2,-r,-2-r-l}(\Gamma ,\Cal S_{++})$ (in $\Cal
S^{1,0,-3-l}(\Gamma ,\Cal S_{++})$ if $r=1$), in local
trivializations.

{\rm (ii)} It follows that when $\psi $ is a morphism from $E_2$ to 
$E_1$, then there are trace expansions for $r>\frac{n+m'+1}2$:
$$\aligned
\Tr [F \psi  (D_{2,\Pi}R^r_{2,\lambda }-D_{1,\Pi}R^r_{1,\lambda })]
&\sim
\sum_{-n\le k< k_0} \tilde c_{ k}(-\lambda )^{\frac{{m'+1}
-k}2-r}\\
&\qquad+ 
\sum_{k\ge k_0}\bigl({ \tilde c'_{ k}}\log (-\lambda
)+{\tilde c''_{ k}}\bigr)(-\lambda )^{\frac 
{{m'}+1 -k}2-r},\\
\Tr [F \psi (D_{2,\Pi}e^{-t\Delta _{2,B}}-D_{1,\Pi}e^{-t\Delta _{1,B}})]
&\sim
\sum_{-n\le k<k_0 } c_{ k}t^{\frac{k-{m'}-1}2}+ 
\sum_{k\ge k_0}\bigl({  c'_{ k}}\log t+{ c''_{ k}}\bigr)t^{\frac 
{k-{m'}-1}2},\\
\Gamma (s)\Tr [F \psi (D_{2,\Pi}\Delta _{2,B}^{-s}-D_{1,\Pi}\Delta _{1,B}^{-s})]
&\sim
\sum_{-n\le k<k_0 } \frac{c_{ k}}{s+\frac{k-{m'}-1}2}\\
&\qquad+
\sum_{k\ge k_0}\Bigl(\frac{-  c'_{ k}}{(s+{\frac 
{k-{m'}-1}2})^2}+\frac{ c''_{ k}}{s+\frac 
{k-{m'}-1}2}\Bigr),
\endaligned\tag3.33$$
with $k_0$ defined by {\rm (3.23)}. For $k\le l-n$, $\tilde c_k$ and
$c_k$ vanish when $k-m'+n$ is even. 

The coefficients $\tilde c_k$, $\tilde c'_k$, $c_k$
and $c'_k$ are locally determined, the $c_k$ and $c'_k$ being
proportional to $\tilde c_k$, $\tilde c'_k$ by universal factors.
\endproclaim  

\demo{Proof} The
statement on the $\psi $do part of (3.32) is immediate, since it is strongly
polyhomogeneous of degree $-2r-1$. For the s.g.o.\ part,
we start by using the
analysis of (3.16) in the proof of Theorem 3.4, now considering the
21-block as in (3.30).  This shows that the symbol-kernel of
$\widetilde G^{(1)}_{\mu }$ is in $\Cal
S^{1,0,-3-l}(\Gamma ,\Cal S_{++})$, in local
trivializations. Next, we apply the composition rules from Theorem 2.7 to
(3.31), combining the above with the information in Theorem 3.4.
It is found that s.g.o.\ part of the first term in (3.31) has
symbol-kernel in $\Cal S^{1,1-r,-2-r-l}(\Gamma ,\Cal S_{++})$, and the
s.g.o.\ part of the second term has symbol-kernel in $\Cal
S^{2,-r,-2-r-l}(\Gamma ,\Cal S_{++})$. This shows (i). It follows
that 
$$
F\psi (D_{2,\Pi}R_{2,-\mu ^2}^r-D_{1,\Pi}R_{1,-\mu ^2}^r)=(F\psi
(D_2Q^r_{2,-\mu 
^2}-D_1 Q^r_{1,-\mu ^2}))_+
+ F\psi \widetilde G^{(r)}_\mu ,
\tag3.34
$$
where the $\psi $do part has symbol in
$S^{1,0,m'-2r}_{\operatorname{spgh,ut}}(\Gamma )$ and $F\psi \widetilde
G^{(r)}_{\mu }$
has symbol-kernel in \linebreak$\Cal S^{1,1-r,m'-2-r-l}(\Gamma ,\Cal S_{++})
+\Cal S^{2,-r,m'-2-r-l}(\Gamma ,\Cal S_{++})$ (with $m'$ moved to the
first upper index when $F$ is tangential).

Now consider (ii).
The trace expansion of the $\psi $do-part is well-known to be a series
of integer powers of $\mu =(-\lambda )^{\frac12}$, beginning with 
$\tilde c_{-n}\mu ^{m'+1+n-2r}$, the terms vanishing when
$k-m'-1+n$ is odd, i.e., $k-m'+n$ is even. 
For the s.g.o.-part, we apply Theorem 2.10.
We get a sum of two $\psi $do terms with $(m,d,s)$ equal
to $(1,1-r,m'-1-r-l)$ resp\. $(2,-r,m'-1-r-l)$ in general (the $m'$
can be moved to the first upper index if $F$ is tangential). They both
give expansions starting with the pure power $\mu ^{m'+n-2r-l}$,
whereas the logarithmic terms start with $\mu ^{m'-2r-l}\log \mu $
resp\. $\mu ^{m'-1-2r-l}\log \mu $; in the case where $F$ is
tangential, the logs start with $\mu ^{-2r-l}\log \mu $
resp\. $\mu ^{-1-2r-l}\log \mu $. Adding the contributions, we find
the first expansion in (3.33).

This carries over to the other two expansions as in Corollary 3.7,
when we furthermore note that $D_i\Pi _0(\Delta _{i,B})=0$ for $i=1,2$.
\qed
\enddemo 

\example {Remark 3.12} Note in particular that $\tilde c_{-n}$ and
$c_{-n}$ in (3.33) vanish when $F=I$.
\endexample

An elaboration of the proof Theorem 3.3 with $F $ replaced by $F\psi D$ gives:
$$
\Tr (F \psi DR^r_\lambda )
\sim
\sum_{-n\le k<0} \tilde b_{ k}(-\lambda ) ^{\frac{{m'}+1 -k}2-r}+ 
\sum_{k\ge 0}\bigl({ \tilde b'_{ k}}\log (-\lambda ) +{\tilde b''_{
k}}\bigr)(-\lambda ) ^{\frac{{m'}+1 -k}2-r},\tag3.35
$$
with $b_{-n}=0$ if $m'$ is even, and with the summation limit $0$
replaced by $m'$ if $F$ is tangential. (When $F$ is tangential,
$m'$ is added to the first upper index instead of the third upper index of the
s.g.o.\ symbol-kernel space.) 
Theorem 3.11 now implies the perturbation result:

\proclaim{Corollary 3.13} Hypotheses and definitions as in Theorem
{\rm 3.11}. For the 
trace expansion {\rm (3.35)} of $F\psi D_1(\Delta _{1,B}-\lambda
)^{-r}$, 
the
 replacement of $D_1$ by
$D_2$ leaves the log-coefficients 
$\tilde b'_k$ invariant for $ k<k_0$. The other coefficients with
$k<k_0$ are modified only by local terms; those with
$ k\le l-n$ and $k-m'+n$ even are invariant.

There are similar results for the associated heat trace and eta function.
\endproclaim

Let us observe a particular consequence for eta expansions. In the above
notation, the eta expansion proved in [GS95], [G99] has the form
$$
\Gamma (s)\Tr [ \psi D_{\Pi}\Delta _{B}^{-s}]
\sim
\sum_{-n< k<0 } \frac{b_{ k}}{s+\frac{k-1}2}+
\sum_{k\ge 0}\Bigl(\frac{-  b'_{ k}}{(s+{\frac 
{k-1}2})^2}+\frac{ b''_{ k}}{s+\frac 
{k-1}2}\Bigr).\tag3.36
$$
With the customary definition $\eta (\psi ,D_\Pi,s')=\Tr(\psi D_{\Pi}\Delta
_B^{\frac{s'+1}2})$, this may also be written in the more well-known form:
$$
\Gamma (\tfrac{s'+1}2)\eta (\psi ,D_\Pi,s')
\sim
\sum_{-n< k<0 } \frac{2b_{ k}}{s'+k}+
\sum_{k\ge 0}\Bigl(\frac{- 4b'_{ k}}{(s'+k)^2}+\frac{ 2b''_{ k}}{s'+k}\Bigr).\tag3.37
$$
We have in a similar way as in Theorem 3.9:

\proclaim{Theorem 3.14} The coefficient $b'_1$ in {\rm (3.36),
(3.37)} is
the same for a non-product type operator $D$
 {\rm (1.1), (1.2)} over $X_c$ with volume form $v(x)\,dx$ and the
associated
product type operator $D^0$
{\rm (1.3)} with volume form $v(x',0)\,dx$, when $P_0=0$ and
$\partial _{x_n}v(x',0)=0$.
Moreover,
the coefficient $b''_1$ differs in the
cases of $D$ and $D^0$
by a local contribution only.
\endproclaim 

There are similar statements for the associated resolvent and heat
trace expansions, as well as extensions to cases where $D-D^0$
vanishes to a general order on $X'$.

In the case $P_0\ne 0$, the coefficient $b'_0$ is invariant under the
replacement of $D$ by $D^0$, and
$b''_0$ is changed only by local contributions; this was shown in
[GS95] for $\Pi =\Pi _>$ plus certain finite rank projections.

For the general systems $\{P-\lambda ,S\varrho \}$ considered in
[G99], one can study the effect of a perturbation of $P$ by a tangential
operator $x_n^l\overline P$ in a similar way, finding also here that
the first $l+1$ logarithmic coefficients are stable, the power coefficients
behind them being changed only by local contributions.

\subhead 4. The resolvent structure for perturbations commuting with $A$ 
\endsubhead

We here consider the case where $D$ is a perturbation of $D^0$ such
that $D-D^0$ {\it commutes with $A$ on} $X_c$,
in the sense 
that in (1.2), the zero-order $x_n$-independent operator (morphism) $P_0$
commutes with $A$, and in the Taylor expansions on $X_c$,$$
x_nP_1(x_n)=\sum_{1\le k\le
K}x_n^kP_{1k}+x_n^{K+1}P'_{K+1}(x_n) \text{ for any }K, 
\tag 4.1$$
the tangential $x_n$-independent first-order differential operators
$P_{1k}$ commute with $A$. The product measure is used on $X_c$.  
We shall show that in this case, there are no log-terms in the
trace expansions in the odd-dimensional case. 

It is no restriction to replace $X_c$ by $X_1$; this can always be
obtained by a scaling in $x_n$.

We know from Theorem 3.5 that the larger $K$ is, the more log-terms are
unaffected by subtracting $x_n^{K+1}P'_{K+1}(x_n)$ from
$x_nP_1(x_n)$, so we may 
disregard this remainder term in the calculations that follow.

Thus, consider the case where, on $X_1$,$$\aligned
D&=\sigma (\partial _{x_n}+A_1(x_n)),\quad A_1(x_n)=A+\sum_{0\le k\le
K}x_n^kP_{1k}=A+\overline P,\\
&\text{ where the $x_n$-independent operators }
P_{1k}\text{ commute with }A.\endaligned
\tag4.2
$$ 
For notational convenience, $P_0$ is here denoted $P_{10}$; it is of
order 0 and the $P_{1k}$ with $k\ge 1$ are of order 1. 

We here restrict the attention to the boundary condition (1.6). 

For the doubled-up systems (cf\. (1.10), (1.16)), we have on
$X_1$, for any $\mu $,$$\aligned 
\Cal D+\mu &=\pmatrix 1&0\\0&\sigma  \endpmatrix \pmatrix \mu & \partial
_{x_n}-A_1^*\\  \partial _{x_n}+A_1 &\mu 
\endpmatrix 
\pmatrix 1&0\\0&\sigma ^* \endpmatrix ,\\
\Cal D^0+\mu  &=\pmatrix 1&0\\0&\sigma  \endpmatrix \pmatrix \mu & \partial
_{x_n}-A\\  \partial _{x_n}+A &\mu
\endpmatrix 
\pmatrix 1&0\\0&\sigma ^* \endpmatrix ,
\endaligned
\tag4.3
$$
since $\left(\smallmatrix 1&0\\ 0&\sigma  
\endsmallmatrix\right)^{-1}=\left(\smallmatrix 1&0\\ 0&\sigma^*  
\endsmallmatrix\right)$. By composition with $\left(\smallmatrix
1&0\\ 0&\sigma   \endsmallmatrix\right)$ and its inverse, the study
of the resolvents 
on $X_1$ is reduced to the study of the inverses of the middle factors
in (4.3), i.e., the case where $\sigma $ is the identity (in $E'_1$).
To keep the notation simple, we use the names $\Cal D+\mu $ and
$\Cal D^0+\mu $ again for the middle factors. In other words,
without of loss of generality: 

{\it We consider in the following the reduced case where $\sigma $ is
the identity, i.e.,$$
\Cal D= \pmatrix 0& \partial
_{x_n}-A_1^*\\  \partial _{x_n}+A_1 &0
\endpmatrix ,\quad
\Cal D^0= \pmatrix 0& \partial
_{x_n}-A\\  \partial _{x_n}+A &0
\endpmatrix.
\tag4.4
$$
The boundary condition {\rm (1.11)} then has the form:}
$$
\Cal B\gamma _0u=0,\quad \Cal B=\pmatrix \Pi _\ge & \Pi _< 
\endpmatrix .\tag4.5
$$ 

Let us denote$$
\Cal P_k=\pmatrix 0&-P_{1k}^*\\ P_{1k}& 0 
\endpmatrix \text{ for }0\le k\le K,\quad \Cal P=\pmatrix
0&-\overline P^*\\ \overline P&0
\endpmatrix, \tag4.6 
$$
then $$
\Cal D-\Cal D^0=\sum_{0\le k\le K}x_n^k\Cal P_k=\Cal P.\tag4.7
$$

We shall use $\zeta _\varepsilon $ introduced in (2.20)ff\. for
$\varepsilon \in \,]0,1]$; it is
defined on 
$X_1$ as constant in $x'$ and
extends by
zero to $X^0$ as well as to $X$, as a $C^\infty $ function.

Rather than $\Cal D=\Cal D^0+\Cal P$, we shall consider $$
\Cal D'=\Cal D^0+\zeta _\varepsilon \Cal P,\tag4.8
$$ 
with $\varepsilon $ to be chosen later; it equals $\Cal D$ on
$X_{\varepsilon /3}$ and serves the same purpose as
$\Cal D$ for investigation of the structure near $x_n=0$.

For $\mu \in \Bbb C\setminus i\Bbb R$, $\Cal D^0 +\mu $ has the
inverse on $\widetilde X^0=X'\times \Bbb R$:$$\aligned
\Cal Q^0 &=(\Cal D^0 +\mu )^{-1}\\
&=\pmatrix \mu (D_{x_n}^2+A^2+\mu ^2)^{-1}&(-\partial
_{x_n}+A)(D_{x_n}^2+A^2+\mu ^2)^{-1}\\ - (\partial
_{x_n}+A)(D_{x_n}^2+A^2+\mu ^2)^{-1}&\mu (D_{x_n}^2+A^2+\mu ^2)^{-1} 
\endpmatrix ,\endaligned\tag4.9
$$
where $(D_{x_n}^2+A^2+\mu ^2)^{-1}=Q^0_\lambda $, $\lambda ={-\mu
^2}$, cf\. (1.17). (The parameter-dependence will not always be
explicitly indicated by an index.) Let$$
\Cal A =\pmatrix \Cal D' +\mu \\ \Cal B \gamma _0
\endpmatrix ,\quad \Cal A^0 =\pmatrix \Cal D^0 +\mu \\ \Cal B \gamma _0 
\endpmatrix,\tag 4.10
$$
representing the full
nonhomogeneous problems on $X^0=X'\times \crp$.
It follows from [GS95, (3.11)--(3.16), Prop\. 3.5] that 
 $\Cal A^0 $ has the solution operator (recall (1.18)ff.) $$
\aligned
(\Cal A^0  )^{-1}&=\pmatrix \Cal R^0 & \Cal K^0  
\endpmatrix,\text{ where }\Cal R^0 =\Cal Q^0_++\Cal G^0 ;\\
\Cal K^0 &=\pmatrix K_{A_\lambda }&0\\0& K_{A_\lambda } 
\endpmatrix S'_{\Cal B},\quad 
S'_{\Cal B}=\pmatrix
\Pi _\ge+\mu^{-1}(A_\lambda+A)\Pi _<\\
\mu^{-1}(A_\lambda-A)\Pi _\ge+\Pi _<\endpmatrix, 
\\ \Cal G^0 &=-\Cal
K^0 \Cal B \gamma _0\Cal Q^0_+=-\pmatrix K_{A_\lambda }&0\\0& K_{A_\lambda } 
\endpmatrix S_{\Cal B} \gamma _0\Cal Q^0_+,\quad S_{\Cal B}=S'_{\Cal
B}\Cal B, 
\endgathered\tag4.11
$$  
for $\mu \in \Bbb C\setminus i\Bbb R$, $\lambda =-\mu ^2$. Here $\Cal
R^0=(\Cal D^0_{\Cal B}+\mu )^{-1}$.
In details, $$\aligned
S_{\Cal B}=S'_{\Cal B}\Cal B &=\pmatrix
\Pi _\ge+\mu^{-1}(A_\lambda+A)\Pi _<\\
\mu^{-1}(A_\lambda-A)\Pi _\ge+\Pi _<\endpmatrix
\pmatrix
\Pi _\ge& \Pi _<\endpmatrix \\
&=\pmatrix
\Pi _\ge&\mu^{-1}(A_\lambda+A)\Pi _<\\
\mu^{-1}(A_\lambda-A)\Pi _\ge&\Pi _<\endpmatrix.
\endaligned\tag4.12
$$
For the description of $\gamma _0\Cal Q^0_+$, we observe the simple
formulas, valid when  $\operatorname{Re}\frak a >0$:
$$\aligned
\tfrac1{\frak a ^2+\xi _n^2}&=\tfrac 1{2\frak a
}(\tfrac1{\frak a +i\xi _n}+\tfrac1{\frak a -i\xi _n}),\text{ so }
h^+\tfrac1{\frak a ^2+\xi _n^2}=\tfrac 1{2\frak a
(\frak a +i\xi _n)},\;h^-\tfrac1{\frak a ^2+\xi
_n^2}=\tfrac1{2\frak a (\frak a -i\xi _n)},\\
\tfrac{i\xi _n}{\frak a ^2+\xi _n^2}&=\tfrac
1{2}(-\tfrac1{\frak a +i\xi _n}+\tfrac1{\frak a -i\xi _n}),\text{ so } 
h^+\tfrac{i\xi _n}{\frak a ^2+\xi _n^2}=-\tfrac 1{2
(\frak a +i\xi _n)},\;h^-\tfrac{i\xi _n}{\frak a ^2+\xi
_n^2}=\tfrac1{2 (\frak a -i\xi _n)};\\
\endaligned\tag 4.13$$
here $h^\pm$ (cf\. (2.6)) projects the rational function
onto its component with poles in 
$\Bbb C_\pm$, respectively. Applying (4.13)
in each eigenspace of $A_\lambda $, we get
(in view of the rules of calculus, cf\. e.g\. [G96, Th\. 2.6.1]):$$
\aligned
\gamma _0(D_{x_n}^2+A^2+\mu
^2)^{-1}_+&=\operatorname{OPT}_n\bigl(h^-\tfrac{1}{A_\lambda ^2+\xi
_n^2}\bigr)=\operatorname{OPT}_n\bigl(\tfrac{1}{2A_\lambda (A_\lambda
-i\xi _n)}\bigr)
=\tfrac1{2A_\lambda }T_{A_\lambda
},\\
 \gamma _0\partial _{x_n}(D_{x_n}^2+A^2+\mu ^2)^{-1}_+
&=\operatorname{OPT}_n\bigl(h^-\tfrac{i\xi _n}{A_\lambda ^2+\xi
_n^2}\bigr)=\operatorname{OPT}_n\bigl(\tfrac{1}{2 (A_\lambda
-i\xi _n)}\bigr)=\tfrac12
T_{A_\lambda },
\endaligned\tag4.14
$$
so that$$
\gamma _0\Cal Q^0_+=\tfrac1{2A_\lambda }\pmatrix {\mu }&
-A_\lambda +A\\ -A_\lambda - A& \mu
\endpmatrix \pmatrix T_{A_\lambda }&0\\0& T_{A_\lambda } 
\endpmatrix =\Cal S^-_1\pmatrix T_{A_\lambda }&0\\0& T_{A_\lambda } 
\endpmatrix ,\tag4.15
$$
with $$
\Cal S_1^-=\tfrac1{2A_\lambda }\pmatrix
{\mu }&
- A_\lambda +A\\ - A_\lambda - A& \mu 
\endpmatrix. \tag4.16
$$
Thus we find from (4.11):$$
\Cal G^0 =\pmatrix K_{A_\lambda }&0\\0& K_{A_\lambda } 
\endpmatrix \Cal S_0 \pmatrix T_{A_\lambda }&0\\0& T_{A_\lambda } 
\endpmatrix ,\quad \text{ with }\Cal S_0=-S_{\Cal B} \Cal
S^-_1.\tag4.17 
$$
In details,$$
\aligned
\Cal S_0&=\tfrac{-1}{2A_\lambda }\pmatrix
\Pi _\ge&\mu^{-1}(A_\lambda+A)\Pi _<\\
\mu^{-1}(A_\lambda-A)\Pi _\ge&\Pi _<\endpmatrix \pmatrix \mu &
-A_\lambda + A\\ -A_\lambda - A& \mu
\endpmatrix\\
&=\tfrac{-1}{2A_\lambda }\pmatrix
\mu \Pi _\ge-\mu ^{-1}(A_\lambda +A) ^2\Pi _<&(-A_\lambda +A)\Pi
_\ge+(A_\lambda+A)\Pi _< \\
(A_\lambda -A)\Pi _\ge-(A_\lambda+A)\Pi _<&-\mu^{-1}(A_\lambda -A)^2\Pi _\ge+\mu \Pi _<\endpmatrix ;
\endaligned\tag4.18
$$
this may be further rewritten by use of (1.5) and the formulas $$
A_\lambda ^2=A^2-\lambda ,\quad (A_\lambda \pm
A)^2=2A^2-\lambda \pm 2A_\lambda A.\tag4.19
$$ 

As shown in Proposition 2.11, $K_{A_\lambda }$ and $T_{A_\lambda }$
are strongly polyhomogeneous Poisson resp\. class 0 trace
operators of degree $-1$ (having symbol-kernels in $\Cal
S^{0,0,-1}(\Gamma ,\Cal S_+)$ in local trivializations). Clearly, $\Cal
S^-_1$ is strongly 
polyhomogeneous of degree $0$ (hence with symbol in $S^{0,0,0}$), so
$\gamma _0\Cal Q^0_+$ is a strongly 
polyhomogeneous trace operator of class 0 and degree $-1$ like
$T_{A_\lambda }$. Since$$
\aligned
\mu ^{-1}(A_\lambda -A)\Pi _{\ge}&=\mu ^{-1}(A_\lambda
-|A|)\tfrac{A_\lambda +|A|}{A_\lambda +|A|}\Pi _{\ge}=\mu (A_\lambda
+|A|)^{-1}\Pi _\ge,\\
\mu ^{-1}(A_\lambda +A)\Pi _{<}&=\mu ^{-1}(A_\lambda
-|A|)\tfrac{A_\lambda +|A|}{A_\lambda +|A|}\Pi _{<}=\mu (A_\lambda
+|A|)^{-1}\Pi _<,
\endaligned\tag4.20
$$ 
where $\Pi _\ge$ and $\Pi _<$ have symbols in $S^0\subset
S^{0,0,0}(\Gamma )$,
$S_{\Cal B} $ is a weakly polyhomogeneous $\psi $do with
symbol in $S^{0,0,0}(\Gamma )$ in local trivializations by Proposition
2.11 (iii) and the product rule in Theorem 2.7 (xiii); then so is
$\Cal S_0$. Thus $\Cal G^0 $ 
is an s.g.o\. with symbol-kernel in $\Cal S^{0,0,-2}(\Gamma ,\Cal
S_{++})$ in local trivializations, in fact with estimates that are
uniform in $x_n\in\Bbb R_+$.

In view of Lemma 2.9 it is,
for the s.g.o.-terms,
the structure near $X'$ that determines their contribution to the
asymptotic expansions; we need not
spend efforts on elaborate presentations of a calculus on the
unbounded manifold $X^0$.

Since $\Cal D^0$ and $\Cal D^0_{\Cal B}$ are skew-selfadjoint (as
unbounded operators in $L_2(\widetilde E^0)$ resp\. $L_2(E^0)$),$$
\|(\operatorname{Re}\mu )\Cal Q^0 \|_{\Cal L(L_2(\widetilde
E^0))}\le C,\;\|(\operatorname{Re}\mu )\Cal R^0 \|_{\Cal
L(L_2(E^0))}\le C,
\text{ for }|\operatorname{Re}\mu |
\ge 1.\tag4.21
$$
In view of the ellipticity, we also have$$
\|\Cal Q^0 \|_{\Cal L(L_2(\widetilde
E^0),H^1(\widetilde
E^0))}\le C',\;\|\Cal R^0 \|_{\Cal L(L_2(
E^0),H^1(
E^0))}\le C',
\text{ for }|\operatorname{Re}\mu |
\ge 1,\tag4.22
$$
with $H^1$ denoting the Sobolev space of order 1.

To find inverses of $\Cal D' +\mu $ and $\Cal A $, we calculate:
$$
\aligned
(\Cal D' +\mu )\Cal Q ^0&=I+\zeta _\varepsilon \Cal P\Cal Q^0 ,
\\
\Cal A (\Cal A ^0)^{-1}&=(\Cal A^0 +\pmatrix \zeta _\varepsilon \Cal P 
\\ 0\endpmatrix )\pmatrix \Cal R^0 &\Cal K^0  
\endpmatrix =I+\pmatrix \zeta _\varepsilon \Cal P\Cal R^0 &\zeta
_\varepsilon \Cal P\Cal K^0 \\0&0 
\endpmatrix 
.\endaligned\tag4.23
$$
Then 
$$\aligned
\Cal Q_{M}&=\Cal Q ^0\sum_{0\le m\le M}(-\zeta _\varepsilon \Cal
P\Cal Q^0 )^m,\\
\Cal C_{M}&=(\Cal A ^0)^{-1}\sum_{0\le m\le M}\pmatrix -\zeta 
_\varepsilon \Cal P\Cal
R^0 &
-\zeta _\varepsilon \Cal P\Cal K^0 \\0&0 
\endpmatrix ^m\\
&= 
\pmatrix \Cal R ^0 &\Cal K^0 
\endpmatrix  \pmatrix \sum_{0\le m\le M}(-\zeta _\varepsilon \Cal
P\Cal R^0 )^m&\sum_{1\le m\le M}(-\zeta _\varepsilon \Cal P\Cal
R^0 )^{m-1}(-\zeta _\varepsilon \Cal P\Cal K^0 ) 
\\0&0 
\endpmatrix \\
&=\Cal R^0 \pmatrix \sum_{0\le m\le M}(-\zeta _\varepsilon \Cal
P\Cal R^0 )^m &-\sum_{0\le m\le M}(-\zeta _\varepsilon \Cal P\Cal
R^0 )^{m-1} \zeta _\varepsilon \Cal P\Cal K^0  
\endpmatrix, \endaligned\tag4.24
$$
will for large $M$ be good approximations to inverses of $\Cal
D' +\mu $ resp\. $\Cal
A $; in particular,
$$
\Cal R_{M}= \Cal R ^0 \sum_{0\le m\le M}(-\zeta _\varepsilon
\Cal P\Cal R^0 )^m \tag4.25
$$
will be a good approximation to
a resolvent of the realization $\Cal D'_{\Cal B }$ of 
$\Cal D' $ under the boundary condition (4.5).
More precisely, we have (cf\. (4.6))
$$
\zeta _\varepsilon \Cal P\Cal R^0 =\zeta _\varepsilon \Cal P_0 \Cal
R^0 +\zeta _\varepsilon x_n\sum_{1\le k\le K}x_n^{k-1}\Cal
P_k\Cal R^0 ,\tag4.26
$$
where the $L_2$ operator norms satisfy (in view of (4.21)--(4.22)): $$
\aligned
\|\Cal P_0\Cal Q^0 \|_{\Cal L(L_2(\widetilde E^0))}\text{ and }\|\Cal P_0\Cal R^0 \|_{\Cal L(L_2(E^0))}&\le
C_1|\operatorname{Re}\mu |^{-1},\\
\|\Cal P_k\Cal Q^0 \|_{\Cal
L(L_2(\widetilde E^0))}\text{ and } \|\Cal P_k\Cal R^0 \|_{\Cal
L(L_2(E^0))}&\le 
C_2,
\endaligned\tag4.27 
$$
for $|\operatorname{Re}\mu |\ge 1$. Since $|\zeta _\varepsilon
x_n|\le \varepsilon $, we can choose an $\varepsilon $ and a $b>0$
such that for 
$|\operatorname{Re}\mu |\ge b$,$$
\|\zeta _\varepsilon \Cal P\Cal Q^0 \|_{\Cal L(L_2(\widetilde
E^0))}\text{ and }
\|\zeta _\varepsilon \Cal P\Cal R^0 \|_{\Cal L(L_2(E^0))}\le
\tfrac12,\tag4.28 
$$
and then there is in fact convergence  in operator norm for
$M\to\infty $:
$$\aligned
\Cal Q &=\lim_{M\to\infty }\Cal Q_{M}=\Cal Q ^0\sum_{ m\ge 0}(-\zeta _\varepsilon \Cal
P\Cal Q^0 )^m,\\
\Cal C &=\lim_{M\to\infty }\Cal C_{M}
=\Cal R^0 \pmatrix \sum_{ m\ge
0}(-\zeta _\varepsilon \Cal 
P\Cal R^0 )^m &-\sum_{ m\ge 1}(-\zeta _\varepsilon \Cal P\Cal
R^0 )^{m-1} \zeta _\varepsilon \Cal P\Cal K^0  
\endpmatrix \\
&=\pmatrix \Cal R  &\Cal K  
\endpmatrix ;\text{ in particular,}\\
\Cal R &=\lim_{M\to\infty }\Cal R_{M}=\Cal R^0  \sum_{ m\ge 0}(-\zeta _\varepsilon \Cal
P\Cal R^0 )^m ;\quad \Cal K =-\Cal R \zeta _\varepsilon \Cal P\Cal K^0  .
\endaligned\tag4.29$$
Here$$
\Cal Q =(\Cal D' +\mu )^{-1},\quad
\Cal C =\Cal A ^{-1},\quad \Cal R =(\Cal D'_{\Cal B }+\mu
)^{-1}=\Cal Q_++\Cal G,\tag4.30 
$$
where $\Cal G$ is an s.g.o\. on $X^0$.
It is seen as in [GS95], [G99], that the operators
belong to the weakly polyhomogeneous calculus (in fact with
estimates that are uniform in $x_n\in \rp$). In particular (as in the
description of $\Cal R_\mu $ in (3.1)ff.), $\Cal G$ has
symbol-kernel in $\Cal S^{0,0,-2}(\Gamma ,\Cal S_{++})$ and $\Cal Q $
has symbol in $ S^{0,0,-1}_{\operatorname{sphg,ut}}(\Gamma )$, in
local trivializations.

At first sight, since $\zeta _\varepsilon \Cal P$ is of order 1, the
terms $\Cal R^0 (\zeta _\varepsilon \Cal P\Cal R^0 )^m$ are all of
order $-1$, so it seems unpractical to use
the series in $m\in \Bbb N$ to get trace expansions. But a closer
inspection shows that only the pseudodifferential
part of each term remains of order $-1$; for the singular Green part,
the order decreases with increasing $m$. 

\proclaim{Proposition 4.1}
{\rm (a)} For each $m$,
$$\aligned
\Cal R^0 (\zeta _\varepsilon \Cal P\Cal R^0 )^m&=(\Cal
Q^0_++\Cal G^0 )\zeta _\varepsilon \Cal P (\Cal Q^0_++\Cal
G^0 )\cdots \zeta _\varepsilon \Cal P (\Cal Q^0_++\Cal
G^0 )\\
&=(\Cal Q^0 (\zeta _\varepsilon \Cal P\Cal Q ^0)^m)_+ +\Cal G_{(m)},
\endaligned\tag4.31
$$
where $\Cal Q^0 (\zeta _\varepsilon \Cal P\Cal Q ^0)^m$ has symbol
in $S^{0,0,-1}_{\operatorname{sphg,ut}}(\Gamma )$
and $\Cal G_{(m)}$ has symbol-kernel in 
\linebreak$\Cal S^{0,0,-m-2}(\Gamma ,\Cal S_{++})$, in local trivializations.

{\rm (b)} With $\Cal Q_M$ and $\Cal R_M$ defined in {\rm (4.24), (4.25)},
one has for any $M\in\Bbb N$: 
$$\Cal R =\Cal
R_{M}+(\Cal Q -\Cal Q_{M})_++\Cal G'_M,
\tag4.32$$
where  $\Cal G'_M$ has
symbol-kernel in $\Cal S^{0,0,-M-3}(\Gamma ,\Cal S_{++})$ and $\Cal
Q -\Cal Q_{M}$ 
has symbol in \linebreak$ S^{0,0,-1}_{\operatorname{sphg,ut}}(\Gamma )$, in
local trivializations.

\endproclaim 

\demo{Proof}
(a) Since $\Cal P\Cal Q^0 $ is strongly polyhomogeneous of order 0, the
statement on the symbol of the $\psi $do part follows
straightforwardly from the product rules; it is the s.g.o.\ part that
demands some effort. 

Note that
$$
\zeta _\varepsilon \Cal P=\zeta _\varepsilon (x_n)\Cal
P_0+x_n\sum_{1\le k\le K}\zeta _\varepsilon (x_n)x_n^{k-1}\Cal P_k,
\tag4.33$$
the sum of a zero-order term
$\zeta _\varepsilon \Cal P_0$ (independent of $\mu $) and a term
containing the factor $x_n$. As we know from Lemma 2.3, a factor
$x_n$ reduces the order of s.g.o.s, lowering the third upper index by 1. Thus
the s.g.o\.  part of $\zeta _\varepsilon \Cal P\Cal R^0 $ has
symbol-kernel in 
$\Cal S^{0,0,-2}(\Gamma ,\Cal S_{++})$ just like $\Cal G^0 $. 
When we multiply out (4.31), the result then follows immediately by
use of Theorem 2.7 (iv)--(vi) for those products that
do not contain two adjacent factors  $\Cal Q^0_+$ and  $\zeta
_\varepsilon \Cal P\Cal Q^0_+$. 

For the remaining products, we need some extra considerations. As in
Theorem 2.7 (xiv), write$$
\Cal Q^0_+\zeta _\varepsilon \Cal P\Cal Q^0_+=
(\Cal Q^0 \zeta _\varepsilon \Cal P\Cal Q^0 )_+-G^+(\Cal
Q^0 )\zeta _\varepsilon \Cal P\,G^-(\Cal Q^0 ) 
$$
(this makes good sense also when $\Bbb R^n$ is replaced
 by $X'\times \Bbb R$). The s.g.o.s $G^\pm (\Cal Q^0 )$ have
symbol-kernels 
in $\Cal S^{0,0,-2}(\Gamma ,\Cal S_{++})$, and thanks to the
structure of $\zeta _\varepsilon \Cal P$ described in (4.33)ff., the
composition $G_1=G^+(\Cal Q^0 )\zeta _\varepsilon \Cal PG^-(\Cal Q^0 )$ has
symbol-kernel in  $\Cal S^{0,0,-3}(\Gamma ,\Cal S_{++})$.
Next, a repeated composition with $\zeta _\varepsilon \Cal P\Cal Q^0_+$
gives$$
\Cal Q^0_+(\zeta _\varepsilon \Cal P\Cal Q^0_+)^2=
((\Cal Q^0 \zeta _\varepsilon \Cal P\Cal Q^0 )_+-G_1)\zeta
_\varepsilon \Cal P\Cal Q^0_+= 
(\Cal Q^0 \zeta _\varepsilon \Cal P\Cal Q^0 )_+\zeta _\varepsilon \Cal P\Cal Q^0_++G_2,
$$
where $G_2$ has symbol-kernel in $\Cal S^{0,0,-4}(\Gamma ,\Cal
S_{++})$ by Lemma 2.3 and Theorem 2.7 (iv)--(vi). Applying
Theorem 2.7 (xiv) to $(\Cal Q^0 \zeta
_\varepsilon \Cal P\Cal Q^0 )_+\zeta _\varepsilon \Cal P\Cal
Q^0_+$, we find that it is the sum of a $\psi $do term $(\Cal
Q^0 (\zeta _\varepsilon \Cal P\Cal Q^0 )^2)_+$ and an s.g.o.\
term $$
G_3=-G^+(\Cal Q^0 \zeta _\varepsilon \Cal P\Cal Q^0 )\zeta
_\varepsilon \Cal P G^-(\Cal Q^0 ).
$$
Inside $G^+$, we apply the commutator formula
$$
x_n\operatorname{OP}(p)=\operatorname{OP}(p)x_n+\operatorname{OP}(i\partial _{\xi _n}p)\tag4.34
$$
to $x_n\Cal Q^0 $, whereby we get two terms, one having a factor
$x_n$ to the right and one where the third upper index
is lowered one step. The 
composition rules and Lemma 2.3 then give that $G_3$ has
symbol-kernel in $\Cal S^{0,0,-4}(\Gamma ,\Cal
S_{++})$. Clearly, this analysis can be continued inductively to 
show that the s.g.o.\ part of $\Cal Q^0_+(\zeta _\varepsilon \Cal
P\Cal Q^0_+)^l$ has symbol-kernel in
$\Cal S^{0,0,-l-2}(\Gamma ,\Cal
S_{++})$ for any $l$, and when this is combined with the other rules, we
obtain (a) for general $m$.

(b) It is clear from (4.24), (4.29), (4.30)ff.\ that the
$\psi $do part of $\Cal R-\Cal R_M$
equals $(\Cal Q -\Cal Q_{M})_+$ and is of order $-1$. For the
s.g.o.\ part, we shall use 
that in view of (4.29) and (a):$$ 
\aligned
\Cal R -\Cal R_{M}&=\Cal R^0 \sum_{ m> M}(-\zeta _\varepsilon \Cal
P\Cal R^0 )^m  
=\Cal R^0 (-\zeta _\varepsilon \Cal
P\Cal R^0 )^M (-\zeta _\varepsilon \Cal P)\Cal R^0  
\sum_{ m\ge 0}(-\zeta _\varepsilon \Cal
P\Cal R^0 )^m  \\
&=(-1)^{M+1}((\Cal Q^0 (\zeta _\varepsilon \Cal P\Cal Q ^0)^M)_+ +\Cal
G_{(M)})\zeta _\varepsilon \Cal P\Cal R .
\endaligned
$$
Here $\Cal R =\Cal Q_++\Cal G$ as described after (4.30). From
the description of $\Cal G_{(M)}$ in (a) follows in view of
Lemma 2.3 and Theorem 2.7 that $\Cal
G_{(M)}\zeta _\varepsilon \Cal P\Cal R $ has symbol-kernel in \linebreak$\Cal
S^{0,0,-M-3}(\Gamma ,\Cal 
S_{++})$. By use of (4.34), we can write $\Cal Q^0 (\zeta
_\varepsilon \Cal P\Cal Q ^0)^M$ as a sum of strongly
polyhomogeneous terms of order $-M-1+j$ with a factor $x_n^j$ to the
right, $j=0,1,\dots ,M$; then it is seen as in the proof of
(a) that the s.g.o.\ part of the composition
$(\Cal Q^0 (\zeta
_\varepsilon \Cal P\Cal Q ^0)^M)_+\zeta _\varepsilon \Cal P\Cal R $
has symbol-kernel in  $\Cal S^{0,0,-M-3}(\Gamma ,\Cal S_{++})$. 
This completes the proof.
\qed
\enddemo 

Since $\Cal G'_M$ has symbol-kernel in $\Cal S^{0,0,-M-3}(\Gamma
,\Cal S_{++})$, it is trace-class for $M>n-3$ and its trace has an
expansion as in Theorem 2.10 beginning with the power $\mu ^{n-M-3}$.
Thus in 
(4.32), the
second term contributes no logarithms in trace expansions and the third
term contributes $O(\mu ^{n-M-2})$ terms with $M$ as large as we
want, so all information on log-terms can be found from the $\Cal R_M$
(for large $M$), and we only have to study  $\Cal R_M$ in detail.

\proclaim{Proposition 4.2} 
 For each $M$, let 
$$\Cal R^0_M=\sum_{0\le m\le M}
\Cal R^0(- \Cal
P\Cal R^0 )^m ;\text{ let
$\Cal G^0_M=$ the s.g.o.\ part of }\Cal R^0_M.\tag4.35
$$ Then
$$
\Cal R_M=\Cal Q_{M,+} +\Cal G^0_M+\Cal G''_M
 ,\tag4.36
$$
where $\Cal G''_M$ has symbol-kernel in $\Cal S^{-\infty ,-\infty
,-\infty }(\Gamma ,\Cal S_{++})$.
\endproclaim

\demo{Proof} 
Denote $\zeta _\varepsilon =\zeta $, $\zeta _{\varepsilon
/3}=\zeta _0$. For each $m\le M$, write$$\aligned
\zeta _0\Cal R^0&(\Cal P \Cal R^0)^m
-\zeta _0\Cal R^0(\zeta \Cal P \Cal
R^0)^m\\
&= \zeta _0\Cal R^0(\Cal P \Cal R^0)^m
-\zeta _0\Cal R^0(\zeta  \Cal P \Cal
R^0)(\Cal P \Cal R^0)^{m-1}+\zeta _0\Cal R^0(\zeta  \Cal P \Cal
R^0)(\Cal P \Cal R^0)^{m-1}\\
&\quad -\zeta _0\Cal R^0(\zeta  \Cal P \Cal
R^0)^2(\Cal P \Cal R^0)^{m-2}+
\zeta _0\Cal R^0(\zeta \Cal P \Cal
R^0)^2(\Cal P \Cal R^0)^{m-2} \\
&\quad -\dots -\zeta _0\Cal R^0(\zeta  \Cal P \Cal
R^0)^{m-1}(\Cal P \Cal R^0)+
\zeta _0\Cal R^0(\zeta  \Cal P \Cal
R^0)^{m-1}(\Cal P \Cal R^0)- \zeta _0\Cal R^0(\zeta  \Cal P \Cal
R^0)^{m}\\
&=\sum_{0\le j\le m-1}\zeta _0\Cal R^0(\zeta \Cal
P\Cal R^0)^j(1-\zeta )\Cal
P\Cal R^0(\Cal
P\Cal R^0)^{m-1-j}.
\endaligned$$
Each term in the sum over $j$ has a factor of the form $\zeta _0(P_+
+ G)(1-\zeta
)$; here $\zeta _0\Cal
P_+(1-\zeta)$ is a $\psi $do with symbol in $S^{-\infty ,-\infty
,-\infty }(\Gamma )$ since $\zeta _0(1-\zeta )=0$, and    $\zeta
_0 G(1-\zeta 
)$ is an s.g.o.\ with symbol-kernel in $\Cal S^{-\infty ,-\infty
,-\infty }(\Gamma , \Cal S_{++})$ by Lemma 2.9; let us call such
operators {\it 
negligible}. It then follows from the composition rules that
$\zeta _0\Cal R^0(\Cal P \Cal R^0)^m
-\zeta _0\Cal R^0(\zeta \Cal P \Cal
R^0)^m$ is negligible, so we get by summation over $m$ that $\zeta
_0\Cal R^0_M-\zeta _0\Cal R_M$ is negligible. By Lemma 2.9, the
s.g.o.\ part of $(1-\zeta _0)\Cal R_M$ is likewise negligible, and so
is $(1-\zeta _0)\Cal G^0_M$. This
implies for the s.g.o.\ parts: $$
\aligned
[\Cal R_M]_{\operatorname{s.g.o.}}&=[\zeta _0\Cal
R_M]_{\operatorname{s.g.o.}}+
[(1-\zeta _0)\Cal R_M]_{\operatorname{s.g.o.}}\\
&=[\zeta _0\Cal
R^0_M]_{\operatorname{s.g.o.}}+\text{ negl.\ s.g.o.}=\zeta _0\Cal G^0_M
+\text{ negl.\ s.g.o.}=\Cal G^0_M
+\text{ negl.\ s.g.o.},
\endaligned
$$
as was to be shown.  
\qed
\enddemo

The proposition shows that the s.g.o.\ part of $\Cal R_M$ equals $\Cal G^0_M$ modulo negligible
terms, so it remains to analyze $\Cal G^0_M$.

We now want to use that the operators $P_{1k}$ commute with the
selfadjoint operator $A$. Then they and their adjoints also commute
with the various 
functions of $A$ appearing in the formulas, such as $A_\lambda $,
$|A|$, $\Pi _\ge$, $e^{-x_nA_\lambda }$, etc. 

Consider, to begin with, the first composition$$
\aligned
\Cal R^0 \Cal P\Cal R^0 &=(\Cal Q^0_++\Cal G^0 )\Cal P (\Cal
Q^0_++\Cal G^0 )\\ 
&=(\Cal Q^0 \Cal P\Cal Q^0 )_+- G^+(\Cal Q^0 )\Cal P\,G^-(\Cal Q^0 )
+  \Cal Q^0_+\Cal P \Cal G^0  +
\Cal G^0 \Cal P \Cal Q^0_++
\Cal G^0 \Cal P \Cal G^0 
\endaligned\tag4.37
$$
where explicit formulas are manageable to some extent.
In the last three terms we shall use that $\Cal G^0 $ has the form
(4.17), where we can write $$\aligned
\sum _{0\le k\le K}x_n^k\Cal P_k\Cal G^0 &=\sum _{0\le k\le K}x_n^k\pmatrix K_{A_\lambda
}&0\\0& K_{A_\lambda }  
\endpmatrix \Cal P_k\Cal S_0 \pmatrix T_{A_\lambda }&0\\0& T_{A_\lambda } 
\endpmatrix\\
\Cal G^0 \sum_{0\le k\le K} x_n^k\Cal P_k&=\sum_{0\le k\le K}\pmatrix K_{A_\lambda
}&0\\0& K_{A_\lambda }  
\endpmatrix \Cal S_0\Cal P_k \pmatrix T_{A_\lambda }&0\\0& T_{A_\lambda } 
\endpmatrix x_n^k,\\
\endaligned\tag4.38$$
by commuting the blocks in the $\Cal P_k$ with $K_{A_\lambda
}=\operatorname{OPK}_n(e^{-x_nA_\lambda })$ resp\. $T_{A_\lambda
}=$ \linebreak$\operatorname{OPT}_n(e^{-x_nA_\lambda })$. 
For the second term we shall use that $G^{\pm}(\Cal
Q^0 )$ have a somewhat similar structure as $\Cal G^0 $.
Some
elementary calculations are needed:

\proclaim{Lemma 4.3} One has for $k, k'\ge 0$:
$$\aligned
\text{\rm (i)}&\quad x_n^k K_{A_\lambda }=\operatorname{OPK}_n(k!(A_\lambda
+i\xi _n)^{-k-1});\quad  T_{A_\lambda
}x_n^k=\operatorname{OPT}_n(k!(A_\lambda -i\xi _n)^{-k-1}).\\ 
\text{\rm (ii)}&\quad T_{A_\lambda }x_n^k K_{A_\lambda }=k!(2A_\lambda
)^{-k-1}.\\ 
\text{\rm (iii)}&\quad \tr_n (x_n^kK_{A_\lambda }ST_{A_\lambda
}x_n^{k'})
=(k+k')!(2A_\lambda )^{-k-k'-1}S,\text{ if }S\text{ commutes with
}A_\lambda .\\
\text{\rm (iv)}&\quad G^{\pm}(\Cal Q^0 )=\pmatrix K_{A_\lambda }&0\\0& K_{A_\lambda } 
\endpmatrix \Cal S_1^\pm \pmatrix T_{A_\lambda }&0\\0& T_{A_\lambda } 
\endpmatrix , \text{ with }\\
&\quad\Cal S_1^\pm=\tfrac1{2A_\lambda }\pmatrix
{\mu }&
\pm A_\lambda +A\\ \pm A_\lambda - A& \mu 
\endpmatrix. \\
\text{\rm (v)}&\quad \pmatrix T_{A_\lambda }&0\\0& T_{A_\lambda } 
\endpmatrix x_n^k \Cal Q^0_+=\sum_{l=0}^{k} \Cal S_{kl}^- \pmatrix T_{A_\lambda }&0\\0& T_{A_\lambda } 
\endpmatrix x_n^{k-l},\\
\text{\rm (vi)}&\quad \Cal Q^0_+x_n^k \pmatrix K_{A_\lambda }&0\\0&
K_{A_\lambda }  \endpmatrix = \sum_{l=0}^{k}x_n^{k-l}\pmatrix K_{A_\lambda
}&0\\0& K_{A_\lambda }  \endpmatrix \Cal S_{kl}^+,
\endaligned$$
where the $\Cal S^\pm_{kl}$ are $2\times2$-matrices whose entries are linear
combinations of $\mu A_\lambda ^{-1-l}$, $AA_\lambda ^{-1-l}$ and
$A_\lambda ^{-l}$.
\endproclaim  
\demo{Proof} 
Rule (i) follows from (1.18) and the formulas$$
(\pm i\partial _{\xi _n})^k(A_\lambda \pm i\xi _n)^{-1}=k!(A_\lambda
\pm i\xi _n)^{-k-1}\tag4.39
$$
and the fact that $$
\aligned
x_n^k\operatorname{OPK}_n( f(\xi
_n))&=\operatorname{OPK}_n((i\partial _{\xi _n})^k f(\xi _n));\\
\operatorname{OPT}_n( f_1(\xi
_n))x_n^k&=\operatorname{OPT}_n((-i\partial _{\xi _n})^k f_1(\xi
_n)).
\endaligned\tag4.40
$$

Rule (ii) follows from:$$
T_{A_\lambda }x_n^kK_{A_\lambda }=\int_0^\infty e^{-x_nA_\lambda
}x_n^ke^{-x_nA_\lambda }\,dx_n=k!(2A_\lambda )^{-k-1}.
$$

Rule (iii) follows from the calculation$$
\multline
\tr_n(x_n^kK_{A_\lambda }ST_{A_\lambda }x_n^{k'})=\int_0^\infty
x_n^ke^{-x_nA_\lambda }Se^{-x_nA_\lambda 
}x_n^{k'}\,dx_n\\
=\int_0^\infty
x_n^{k+k'}e^{-x_n2A_\lambda }S\,dx_n
=(k+k')!(2A_\lambda )^{-k-k'-1}S.
\endmultline\tag4.41
$$

For (iv), we use that when $p(\xi _n)$ is the symbol of a $\psi $do
$P$ of order $\le 0$, then the symbol-kernel of $G^{\pm}(P)$ equals
$[\Cal F^{-1}_{\xi _n\to z_n}h^\pm p(\xi _n)]_{z_n=\pm (x_n+ y_n)}$;
cf\. e.g\. [G96, Th\. 2.6.10]. Using this in each eigenspace of
$A_\lambda $, we find in view of (4.13): $G^+(\frac1{A_\lambda
^2+D_{x_n}^2})$ has the symbol-kernel$$
\bigl[\Cal F^{-1}_{\xi _n\to z_n}\tfrac1{2A_\lambda (A_\lambda +i\xi
_n)}\bigr]_{z_n=x_n+y_n}=\tfrac 1{2A_\lambda }e^{-(x_n+y_n)A_\lambda };
$$
$G^+(\frac{\partial _{x_n}}{A_\lambda
^2+D_{x_n}^2})$ has the symbol-kernel$$
\bigl[\Cal F^{-1}_{\xi _n\to z_n}\tfrac{-1}{2 (A_\lambda +i\xi
_n)}\bigr]_{z_n=x_n+y_n}=-\tfrac 1{2}e^{-(x_n+y_n)A_\lambda }.
$$
Similarly, $
G^-(\frac1{A_\lambda
^2+D_{x_n}^2})$ has the symbol-kernel
$\tfrac 1{2A_\lambda }e^{-(x_n+y_n)A_\lambda }$ and
$G^-(\frac{\partial _{x_n}}{A_\lambda
^2+D_{x_n}^2})$ has the symbol-kernel
$\tfrac 1{2}e^{-(x_n+y_n)A_\lambda }$. 
In other words,
$$
G^{\pm}(\tfrac1{A_\lambda
^2+D_{x_n}^2})=K_{A_\lambda }\tfrac1{2A_\lambda }T_{A_\lambda },\quad
G^{\pm}(\tfrac{\partial _{x_n}}{A_\lambda
^2+D_{x_n}^2})=\mp\tfrac12 K_{A_\lambda }T_{A_\lambda }.
\tag 4.42$$
Application of these informations to $\Cal Q^0 $ shows (iv).

In the proof of (v) and (vi), we need some further decompositions of
rational functions (cf\. also (4.13)):
$$\aligned
\tfrac1{(\frak a \pm i\xi _n)^k(\frak a \mp i\xi _n)}&=\tfrac1{(\frak a
+i\xi _n)^{k-1}2\frak a }\bigl(\tfrac1{\frak a + i\xi
_n}+\tfrac1{\frak a - i\xi _n}\bigr)
=\tfrac1{2\frak a(\frak a \pm i\xi _n)^k}+\tfrac1{2\frak a
(\frak a \pm i\xi _n)^{k-1}(\frak a \mp i\xi
_n)}\\
&=\dots=\tsize\sum_{j=1}^k\tfrac1{(2\frak a )^{j}(\frak a \pm i\xi _n)^{k+1-j}}
+\tfrac1{(2\frak a )^{k}(\frak a \mp i\xi _n)},\endaligned$$
and hence
$$
\aligned
h^\pm\tfrac 1{(\frak a \pm i\xi _n)^k(\frak a ^2+\xi
_n^2)}&=
h^\pm\tfrac 1{(\frak a \pm i\xi _n)^{k+1}(\frak a \mp i\xi 
_n)}=
\tsize\sum_{j=1}^{k+1}
\tfrac1{(2\frak a )^{j}(\frak a \pm i\xi _n)^{k+2-j}},\\
h^\pm\tfrac {i\xi _n}{(\frak a \pm i\xi _n)^k(\frak a ^2+\xi _n^2)}
&
=h^\pm\Bigl[\tfrac1{(\frak a
\pm i\xi _n)^{k} }\tfrac12\bigl(-\tfrac1{\frak a + i\xi
_n}+\tfrac1{\frak a - i\xi _n}\bigr)\Bigr]\\
&=\mp\tfrac1{2(\frak a \pm i\xi _n)^{k+1}}\pm \tsize\sum_{j=1}^{k}
\tfrac1{2(2\frak a )^{j}(\frak a \pm i\xi _n)^{k+1-j}}.
\endaligned\tag4.43$$
We then get by use of (i) and the rules of calculus (cf\. e.g\.
[G96, Th\. 2.6.1]):
$$
\aligned
T_{A_\lambda }x_n^k(A_\lambda ^2+D_{x_n}^2)^{-1}_+&=
\operatorname{OPT}_n(h^-(\tfrac {k!}{(A_\lambda -i\xi
_n)^{k+1}}\tfrac1{A_\lambda ^2+\xi _n^2}))\\
=
\operatorname{OPT}_n\bigl(
\tsize\sum_{j=1}^{k+1}
&\tfrac{k!}{(2A_\lambda )^{j}(A_\lambda - i\xi _n)^{k+2-j}}\bigr)
=\tsize\sum_{j=1}^{k+1}\tfrac{k!}{(k+1-j)!(2A_\lambda
)^j}T_{A_\lambda }x_n^{k+1-j},\\
(A_\lambda ^2+D_{x_n}^2)^{-1}_+x_n^kK_{A_\lambda }&=
\operatorname{OPK}_n(h^+(\tfrac1{A_\lambda ^2+\xi _n^2}\tfrac
{k!}{(A_\lambda +i\xi 
_n)^{k+1}}))\\
&=\tsize\sum_{j=1}^{k+1}\tfrac{k!}{(k+1-j)!}x_n^{k+1-j}K_{A_\lambda
}\tfrac1{(2A_\lambda )^j},
\endaligned\tag4.44$$
and similarly
$$\aligned
T_{A_\lambda }x_n^k\partial _{x_n}(A_\lambda
^2+D_{x_n}^2)^{-1}_+&
=\operatorname{OPT}_n\bigl(\tfrac1{2(A_\lambda -i\xi _n)^{k+1}}-
\tsize\sum_{j=1}^{k}
\tfrac1{2(2A_\lambda )^{j}(A_\lambda - i\xi _n)^{k+1-j}}\bigr)\\
&=\tfrac {k!}2T_{A_\lambda
}x_n^{k}-\tsize\sum_{j=1}^{k}\tfrac{k!}{2(k-j)!(2A_\lambda
)^j}T_{A_\lambda }x_n^{k-j},
\\
\partial _{x_n}(A_\lambda ^2+D_{x_n}^2)^{-1}_+x_n^kK_{A_\lambda }&=
\operatorname{OPK}_n\bigl(-\tfrac1{2(A_\lambda +i\xi _n)^{k+1}}+
\tsize\sum_{j=1}^{k}
\tfrac1{2(2A_\lambda )^{j}(A_\lambda + i\xi _n)^{k+1-j}}\bigr)\\
&=-\tfrac{ k!}2x_n^{k}K_{A_\lambda
}+\tsize\sum_{j=1}^{k}\tfrac{k!}{2(k-j)!}x_n^{k-j}K_{A_\lambda
}\tfrac1{(2A_\lambda 
)^j}.
\endaligned\tag4.45
$$

Application to the blocks in $\Cal Q^0 $ give formulas (v)
and (vi), with the asserted structure of the $\Cal S^\pm_{kl}$.
\qed

\enddemo 

\example {Remark 4.4}
There are similar results with $\Cal Q^0$ replaced by $(\Cal Q^0)^r$.
The powers contain
factors $(D_{x_n}^2+A^2+\mu ^2)^{-j}$ with higher $j$ which, in the
application of $h^+$ and $h^-$ as in (4.13), (4.43) lead to fractions with
both $A_\lambda + i\xi _n$ and $A_\lambda - i\xi _n$ in higher powers in the denominator. Again one
decomposes into simple fractions so that the numerators are independent
of $\xi _n$, which leads to
formulas generalizing (iv)--(vi) and of a similar form (with $\mu $
and $A$ appearing in higher powers).
\endexample

Using (iv), we can write, similarly to (4.38):$$
G^+(\Cal Q^0 )\sum_{0\le k\le K} x_n^k\Cal P_k=\sum_{0\le k\le K}\pmatrix K_{A_\lambda
}&0\\0& K_{A_\lambda }  
\endpmatrix \Cal S^+_1\Cal P_k \pmatrix T_{A_\lambda }&0\\0& T_{A_\lambda } 
\endpmatrix x_n^k,\tag4.46
$$
For simplicity, {\it we write from now on the diagonal block
matrices formed of 
$K_{A_\lambda }$ or $T_{A_\lambda }$ as simple factors,} meaning
that they are composed with each block (as already done with e.g\.
$1/(2A_\lambda )$); this should not lead to confusion.
The s.g.o\. terms in
(4.33) can now be calculated:

\proclaim{Proposition 4.5} The singular Green part $\Cal G^0_{(1)}$ of
$\Cal R^0 \Cal P\Cal R^0 $ is a 
a block matrix, cf. {\rm (4.52)},
$$
\Cal G^0_{(1)}=\pmatrix \Cal G^0_{(1),11}&  \Cal G^0_{(1),12}\\  \Cal G^0_{(1),21}&
\Cal G^0_{(1),22}  
\endpmatrix ,\tag4.47
$$
where each block is a sum of terms 
$$
x_n^{k'}K_{A_\lambda }S\,T_{A_\lambda }x_n^{k''},\tag4.48
$$ with $S$ of the form$$
\text{\rm (a) }(-\lambda )^{l/2}L_mA_\lambda ^{-j}\tfrac A{|A'|},\quad 
\text{\rm (b) }(-\lambda )^{l/2}L_mA_\lambda ^{-j}\;\text{ or }\;\text{\rm (c) } 
(-\lambda )^{l/2}L_m\Pi _0;
\tag4.49
$$
here $l\in\Bbb Z$, $m$ and $j\in\Bbb N$, and $L_m$ denotes a
($\lambda $-independent) 
differential operator of order $m$ commuting with $A$.
Consequently, the normal trace
is a block matrix, cf\. {\rm (4.54)},$$
\tr_n \Cal G^0_{(1)}=S_{(1)}=\pmatrix S_{(1),11}&  S_{(1),12}\\
S_{(1),21}&  S_{(1),22}   
\endpmatrix ,\tag4.50
$$
where each block 
is a linear combination of 
terms of the form {\rm (4.49)}.
\endproclaim

\demo{Proof} 
We have for the term $G^+(\Cal Q^0 )\Cal PG^-(\Cal Q^0 )$:
$$
\aligned
G^+(\Cal Q^0 )\Cal PG^-(\Cal Q^0 )&=K_{A_\lambda }\Cal S^+_1T_{A_\lambda
}\sum_{0\le k\le K}x_n^k\Cal P_k K_{A_\lambda }\Cal S^-_1
T_{A_\lambda }\\
&=K_{A_\lambda }\sum_{0\le k\le K}k!(2A_\lambda
)^{-k-1}\Cal S^+_1\Cal P_k\Cal S^-_1 
T_{A_\lambda },
\endaligned\tag4.51$$
where we used Lemma 4.3 (iv), (4.46) and Lemma 4.3 (ii). 
Treating
$\Cal G^0 \Cal P\Cal G^0$ in the same way, using (4.38), we get a
similar expression with $\Cal S_0$ instead of $\Cal S_1^\pm$.
For the last two terms $\Cal G^0 \Cal P\Cal Q^0_+$ and $\Cal
Q^0_+\Cal P\Cal G^0$, we moreover use (v) and (vi) in Lemma 4.3,
finding e.g.$$
\Cal G^0 \Cal P\Cal Q^0_+=\sum_{0\le k\le K}K_{A_\lambda }\Cal
S_0\Cal P_kT_{A_\lambda }x_n^k\Cal Q^0_+
=\sum_{0\le k\le K}K_{A_\lambda }\Cal
S_0\Cal P_k\sum_{0\le l\le k}S_{kl}^-T_{A_\lambda }x_n^{k-l}.
$$ 
This gives:$$
\aligned
\Cal G^0_{(1)}&=- G^+(\Cal Q^0 )\Cal P\,G^-(\Cal Q^0 )
+  \Cal Q^0_+\Cal P \Cal G^0  +
\Cal G^0 \Cal P \Cal Q^0_++
\Cal G^0 \Cal P \Cal G^0 \\
&=K_{A_\lambda }\sum_{0\le k\le K}k!(2A_\lambda
)^{-k-1}(-\Cal S^+_1\Cal P_k\Cal S^-_1 +\Cal S_0\Cal P_k\Cal S_0) 
T_{A_\lambda }\\
&\quad+\sum_{0\le k\le K}\sum_{0\le l\le k}[K_{A_\lambda }\Cal
S_0\Cal P_k\Cal S_{kl}^-T_{A_\lambda
}x_n^{k-l}+x_n^{k-l}K_{A_\lambda }\Cal S^+_{kl}\Cal P_k\Cal
S_0T_{A_\lambda }] .  
\endaligned\tag4.52
$$
This shows the general structure (4.48) of the terms in the blocks,
and we shall 
now show the additional information given in (4.49).

The blocks in
$\Cal S^\pm_1$ and the $\Cal S^\pm_{kl}$ are linear combination of terms
as in (4.49)(b).  
$\Cal S_0$ moreover contains terms of the form $$
(-\lambda )^{l/2}L_mA_\lambda ^{-j}\Pi _\ge \tag4.53
$$ (we insert $\Pi _<=I-\Pi _\ge$, use the reductions in (4.19) and
absorb powers of $A$ in $L_m$, noting that since $A_\lambda
=(A^2-\lambda )A_\lambda ^{-1}$, $A_\lambda $ need
only occur explicitly in negative powers). Then in the
resulting matrices, we get 
linear combinations of terms as in (4.49)(b) and (4.53), by moving the
differential operators $P_{1k}$ out in front in each block by
commutation. To the terms of the form (4.53) we apply (1.5), which
leads to$$
(-\lambda )^{l/2}L_mA_\lambda ^{-j}\Pi _\ge =\tfrac12(-\lambda
)^{l/2}L_mA_\lambda ^{-j}(\tfrac A{|A'|}+I+\Pi _0),
$$
giving the three types in (4.49).

(4.50)ff\. follows by application of
Lemma 4.3 (iii) to each block. More precisely, this gives
$$
\aligned
\tr_n \Cal G^0_{(1)}
&=\sum_{0\le k\le K}k!(2A_\lambda
)^{-k-2}(-\Cal S^+_1\Cal P_k\Cal S^-_1 +\Cal S_0\Cal P_k\Cal S_0) 
\\
&\quad+\sum_{0\le k\le K}\sum_{0\le l\le k}(k-l)!(2A_\lambda
)^{-k+l-1}(\Cal S_0\Cal P_k\Cal S_{kl}^-+\Cal S^+_{kl}\Cal P_k\Cal
S_0).\quad\square  
\endaligned\tag4.54
$$
\enddemo 

Note that although $\Cal P$ itself does not commute with $\Cal S_0$,
$\Cal S^\pm_1$, etc.,
 we obtained the result by commutation in each block.

This shows the first step in 

\proclaim{Theorem 4.6} For any $M\ge 0$, the s.g.o\. part $\Cal
G^0_M$ of $\Cal R^0_M $ is of the form
$$
\Cal G^0_{M} =\pmatrix \Cal G^0_{M,11}&  \Cal G^0_{M,12}\\  \Cal G^0_{M,21}&
\Cal G^0_{M,22}   
\endpmatrix, \tag4.55
$$ where the blocks have the structure in 
{\rm (4.48)--(4.49)}.
\endproclaim 

\demo{Proof}
We already have this structure for the s.g.o.\ parts of $\Cal R^0$ and $\Cal
R^0\Cal P\Cal R^0$.
Now consider $\Cal R^0 (\Cal P\Cal R^0 )^m$. Again we depart from
the exact formulas for $\Cal R^0 =\Cal Q^0_++\Cal G^0 $ given
in (4.9) and (4.17)--(4.18), as in the proof of Proposition 4.5, 
using the description in Lemma 4.3 of the effects of
multiplication by $x_n^k$. We moreover use the formulas for higher
powers, as described in Remark 4.4.
Then we find the structure in (4.48)--(4.49) also for the $m$'th
term and the result for $\Cal G^0_M$ follows by summation.
\qed
\enddemo 

Similar results can be shown for powers of the resolvent and for
powers of the blocks in the resolvent. 

{\it We now pass to the consequences for the original
operators} (4.3) {\it with general} $\sigma $, by composing suitably
with $\left(\smallmatrix
1&0\\ 0&\sigma   \endsmallmatrix\right)$ and its inverse. This gives
for the  resolvent
$(\Delta _B-\lambda )^{-1}$:

\proclaim{Theorem 4.7} Under the assumption of {\rm (4.2)}, the
resolvent $R_{\lambda }$ of $\Delta _B=D_\ge^*D_\ge$ satisfies, for
any $M$, any $r\ge 1$:
$$
R^r_\lambda =Q^r_{\lambda ,+}+\zeta _{\varepsilon  }G_{M,r}\zeta
_{\varepsilon  }+G'_{M,r},\tag4.56
$$
where $G_{M,r}$ is a finite sum of terms as in {\rm (4.48)--(4.49)}
and $G'_{M,r}$ has symbol-kernel in 
$\Cal S^{0,-r,-M-r-2}(\Gamma
,\Cal S_{++})$, in local trivializations. 
\endproclaim 

\demo{Proof} As usual, $\lambda =-\mu ^2$. We have on $X_{\varepsilon
 }$, in view
of (3.2) and (4.3), 
$$\multline
R_{-\mu ^2}\sim \mu ^{-1}\pmatrix 1&0  
\endpmatrix \pmatrix 1&0\\ 0&\sigma 
\endpmatrix \Cal R   \pmatrix 1&0\\ 0&\sigma ^* 
\endpmatrix \pmatrix 1\\0
\endpmatrix  \\
=\mu ^{-1}\Cal R_{11}=\mu ^{-1}(\Cal Q_{11,+}+\Cal
G^0_{M,11} +\Cal G'_{M,11}) 
,\endmultline\tag4.57
$$
where $\Cal G^0_{M,11}$ has the structure described in Theorem 4.6 and
$\Cal G'_{M,11}$ has symbol-kernel in $\Cal S^{0,0,-M-3}(\Gamma ,\Cal
S_{++})$. Here $\mu ^{-1}\Cal Q_{11}=Q_{\lambda }$, cf\. (3.2).
This implies the statement for $r=1$ by a couple of applications of
Lemma 2.9; on one hand it
allows the multiplication by $\zeta _{\varepsilon  }$, on the other
hand it allows extending the structure from $X_{\varepsilon  }$ to $X$.

For the higher powers, note that by (4.29), $$
\Cal Q =\Cal Q_{M}+\Cal Q'_{M}, \text{ where }
\Cal Q'_{M}= \Cal Q ^0\sum_{ m>M}(-\zeta _\varepsilon \Cal
P\Cal Q^0 )^m=\Cal Q (-\zeta _\varepsilon \Cal
P\Cal Q^0 )^{M+1}.\tag4.58
$$ 
 Consider the $r$'th power (on $X_{\varepsilon  }$)
$$
R_\lambda ^r\sim [\mu ^{-1}(\Cal Q_{11,+}+\Cal G_{11})]^r= \mu ^{-r}(
\Cal Q_{M,11,+}+\Cal Q'_{M,11,+}+\Cal
G^0_{M,11} +\Cal G'_{M,11})^r.\tag4.59 
$$
The s.g.o.\ part of $R_\lambda ^r$ comes partly from compositions containing
$\Cal G_{11}$,  partly from ``leftover''
contributions from the 
$\psi $do terms (as in Theorem 2.7 (xiv)). The result is obtained by
using again the exact 
formulas given above for the terms in the sums 
over $m$ and the information on the remainders (for $\Cal Q'_{M}$
it is the fact that it contains $M+1$ factors $\Cal P$), in a similar way as
in the preceding proofs.
\qed
\enddemo 

\example{Remark 4.8} It would also have been possible to get this
result --- and even a 
slightly better one with the s.g.o.\ remainder symbol-kernel in
$\Cal S^{0,0,-M-2r-2}(\Gamma
,\Cal S_{++})$ --- by departing from the formulas for $R_\lambda $ in
[G92], but the 
exact rules that would have to be worked out, would be even more
complicated, since the difference between the second order operators
$D^*D$ and ${D^0}'D^0$ 
contains many more terms, including some of the form $\partial _{x_n}P$ and
$x_n\partial _{x_n}P$, and the two second-order realizations have
different 
boundary conditions when $P_{10}\ne 0$.
\endexample 

\example{Remark 4.9} The other resolvent $(D_\ge D^*_\ge-\lambda
)^{-1}$ has a similar form on $X_1$, except that it is composed with
$\sigma $ to the left and $\sigma ^*$ to the right.
\endexample 

\subhead 5. Trace results in the commuting case \endsubhead

In this section, we continue the study of the operator families
defined from $\Delta
_B=D^*_\ge D_\ge$.
For the purpose of analyzing the traces, we introduce a notation for symbols
with alternating parity in the symbol expansion:

\proclaim{Definition 5.1} Let $p(x,\xi )\sim \sum_{l\ge
0}p_{m-l}(x,\xi )$ be a classical $\psi $do symbol of integer order $m$,
expanded in terms $p_{m-l}$ that are homogeneous of degree $m-l$ for
$|\xi |\ge 1$. 

We say that $p$  has {\bf even-even alternating parity}, when the
terms $p_{m-l}$ of even degree $m-l$ are even in $\xi $, and the
terms $p_{m-l}$ of odd degree $m-l$ are odd in $\xi $,
i.e.,$$
p_{m-l}(x,-\xi )=(-1)^{m-l}p_{m-l}(x,\xi ),\text{ for all }l,\tag 5.1
$$
and the same holds for derivatives of $p$.

We say that $p$  has {\bf even-odd alternating parity}, when the
terms $p_{m-l}$ of even degree $m-l$ are odd in $\xi $, and the
terms $p_{m-l}$ of odd degree $m-l$ are even in $\xi $,
i.e.,$$
p_{m-l}(x,-\xi )=(-1)^{m-l+1}p_{m-l}(x,\xi ),\text{ for all }l,\tag 5.2
$$
and the same holds for derivatives of $p$.
\endproclaim 

Note that symbols of differential operators, and parametrix symbols
for elliptic differential operators, have even-even
alternating parity. On the other hand, the symbol of $|A|$ has
even-odd alternating parity; this can be checked on the basis of the 
formulas for the symbol of $(A^2)^{\frac12}$ given in Seeley [S67].
Then also $\frac A{|A'|}$ has even-odd alternating parity.
(Such alternating parity properties were also observed in
[GS96, (3.9)ff.].)

\proclaim{Theorem 5.2}
Let $S$ be a $\mu $-dependent $\psi $do on $X'$ with polyhomogeneous symbol in
$S^{m,d,s}(\Gamma )$ ($s\le 0$) in local trivializations, holomorphic in
$\mu $, and such that 
the homogeneous terms of degree $m+s+d-j$ with $m+s>-n$ are
integrable in $\xi '$. Assume moreover that $S$ is a finite sum of
terms $S_i$, $i=1,\dots,i_3$, of the form
$$
S_i=\cases (-\lambda )^{l_i/2}L_{m_i}A_\lambda ^{-j_i}\tfrac
A{|A'|}\text{ for }i=1,\dots,i_1,\\ 
(-\lambda )^{l_i/2}L_{m_i}A_\lambda ^{-j_i}\text{ for }i=i_1+1,\dots,i_2,\\  
(-\lambda )^{l_i/2}L_{m_i}\Pi _0\text{ for }i=i_2+1,\dots,i_3,\endcases  
\tag5.3
$$
where $l_i\in\Bbb Z$, $m_i$ and $j_i\in\Bbb N$, and $L_{m_i}$ is a
($\lambda $-independent) 
differential operator of order $m_i$.
Then in the asymptotic expansion 
$$
\Tr_{X'}S\sim \sum_{j\in\Bbb N}c_j(-\lambda )^{(m+d+s+n-1-j)/2}+\sum_{k\in\Bbb
N}(c'_k\log(-\lambda ) +c''_k)(-\lambda ) ^{{d+s-k}/2},\tag5.4
$$
the logarithmic terms $c'_k(-\lambda )^t\log(-\lambda )$ with integer
$t$ come from the terms $S_i$ with 
$l_i-j_i$ {\bf even}, $i\le i_1$, and the logarithmic terms with
noninteger $t$ 
($t-\frac12$ integer) come from the terms $S_i$ with
$l_i-j_i$ {\bf odd}, $i\le i_1$. The logarithmic terms
coming from each such $S_i$ have the powers $t= (l_i-j_i)/2-\nu $,
$\nu =0,1,2,\dots$ (when $j_i=0$, there is only the power $l_i/2$).

Moreover, if $n$ is odd, all the
logarithmic terms are zero.
\endproclaim 

\demo{Proof}
Consider first the terms with $i>i_2$, the last type in (5.3). Such
a term has 
smooth finite dimensional range (in particular, it is of order
$-\infty $) and  
contributes a constant $\Tr_{X'}(L_{m_i}\Pi _0)$ times $(-\lambda
)^{l_i/2}$. We start by subtracting these terms from the expression to
be analyzed, which leaves us with a decomposition in terms of the
first two types.

Note that {\it we have not assumed the first two types of terms in
{\rm (5.3)} 
to be 
trace-class operators individually,} and that they need not be turned
into trace-class 
operators  by differentiatiation of high order in $\lambda
$ (because in the Leibniz formula, some differentiations would fall on
$(-\lambda )^{l_i/2}$, others on $A_\lambda ^{-j_i}$).

However, as recalled in the proof of Theorem 2.10,
the regions $\{|\xi '|\ge |\lambda |^\frac12\}$
and $\{|\xi '|\le 1\}$ give 
pure powers; it is the region $\{1\le |\xi
'|\le |\lambda |^{\frac12}\}$ that may contribute with log-power terms.
Only in this region will the decomposition be used; it
corresponds to a similar decomposition 
for the symbols of the operators. The point is now that although the
individual terms here need not be of 
sufficiently low order to allow integration over $\Bbb R^{n-1}$, we can
certainly integrate them over $\{1\le |\xi '|\le |\lambda
|^{\frac12}\}$. 

The terms of with $i_1<i\le i_2$
in (5.3) are of the form of a power of $-\lambda $ times a 
strongly polyhomogeneous operator; their symbols will contribute pure
powers (since they obviously do so when integrated over $\{0\le |\xi '|\le
|\lambda |^{\frac12}\}$, and the region $\{|\xi '|\le 1\}$ gives only
pure powers). 

It is the symbols with $i\le i_1$ in (5.3) that
may
contribute logarithmic terms. We drop the index $i$ in the
following.

Consider such a term $(-\lambda )^{l/2}L_mA_\lambda ^{-j}\tfrac A{|A'|}$.
We may write it:$$
(-\lambda )^{\frac l2}L_mA_\lambda ^{-j}\tfrac A {|A'|}
=\varrho ^{\frac{j-l}2}L_m(\varrho A^2+1)^{-\frac j2}\tfrac A {|A'|},
\quad\varrho 
=-\lambda 
^{-1}. \tag5.5 
$$
If $j>0$, we insert a power series expansion of the $-\frac j2$-power, 
$$
\varrho ^{\frac{j-l}2}L_m(\varrho A^2+1)^{-\frac j2}\tfrac A
{|A'|}\sim \varrho 
^{\frac{j-l}2}
L_m\sum_{\nu \ge 0}\tbinom{-\frac j2}{\nu }\varrho ^{\nu }A^{2\nu
}\tfrac A {|A'|}.\tag5.6 
$$ 
From this formula for the full operators, we can also find the
structure of the symbol by 
inserting the polyhomogeneous symbol
expansions and carry out symbol compositions.
(This type of expansion is somewhat like the Taylor expansion
in [GS95, Th\. 1.12]; in the present case the order of the coefficient of
$\varrho ^\nu $ increases by 2 when $\nu $ increases by 1.
Systematic calculi 
with such features are worked out in Loya [L01] and [GH02].)  

Let $B^{m,\nu }=\tbinom{-\frac j2}{\nu }L_mA^{2\nu }\tfrac A {|A'|}$,
it is of order $m+2\nu $ and its 
symbol $b^{m,\nu }(x',\xi ')$ 
has an expansion$$
b^{m,\nu }(x',\xi ')\sim \sum_{\nu ' \ge 0}b^{m,\nu }_{m+2\nu -\nu '}(x',\xi
'),\tag5.7
$$
where $b^{m,\nu }_{m+2\nu -\nu '}$ is homogeneous of degree $m+2\nu
-\nu ' $ in $\xi '$. It suffices to consider $\lambda \in \Bbb R_-$ (the
results extend analytically to other $\lambda $ in view of [GS95,
Lemma 2.3]). Since $\xi '$ runs in dimension $n-1$,
 $$\aligned
\int_{1\le |\xi '|\le |\lambda |^{\frac12}}&b^{m,\nu }_{m+2\nu -\nu
'} (x',\xi 
')\,\d\xi '\\
&=(2\pi )^{1-n}\int_1^{|\lambda |^{\frac12}}r^{m+2\nu -\nu '+n-2}\,dr
\int_{|\xi '|=1}b^{m,\nu }_{m+2\nu -\nu '} (x',\xi
')\,d\sigma (\xi ')\\
&=\cases c_{\nu '}(x')(|\lambda |^{m+2\nu -\nu '+n-1}-1)&\text{ if
}m+2\nu -\nu 
'+n-1\ne 0,\\
c_{\nu '}(x')\log |\lambda |&\text{ if }m+2\nu -\nu
'+n-1= 0,\endcases
\endaligned\tag5.8
$$
where $$
c_{\nu '}(x')=c \int_{|\xi '|=1}b^{m,\nu }_{m+2\nu -\nu '} (x',\xi
')\,d\sigma (\xi '),\tag5.9
$$
with a nonzero constant $c $ (depending on $m+2\nu
-\nu '+n-1$).
The case of the logarithm occurs when $\nu '=m+2\nu +n-1$, and then
$$
c_{m+2\nu +n-1}(x')=c \int_{|\xi '|=1}b^{m,\nu }_{1-n} (x',\xi
')\,d\sigma (\xi ').\tag5.10
$$ 

We see that each term in the expansion (5.7) contributes with one
logarithmic term, and that it is proportional to  $(-\lambda
)^{\frac {l-j}2-\nu }\log(-\lambda )c_{m+2\nu +n-1}(x')$. The power
is integer resp\. integer$+\frac12$ exactly when $(l-j)/2$ is
integer resp\. integer$+\frac12$, i.e., when $l-j$ is even resp\. odd.
The highest order log-term comes from the case $\nu =0$ and 
is of the form $$
c(x')(-\lambda )^{\frac {l-j}2}\log(-\lambda ).\tag5.11
$$

If $j=0$ in (5.5), we have just one term $(-\lambda )^{\frac
l2}L_m\frac A{|A'|}$ to analyze. Studying $L_m\frac A{|A'|}$ as we
did with $B^{m,\nu }$, we find a single logarithmic contribution
(from the term of degree $1-n$ in the symbol)
$$
c(x')(-\lambda )^{\frac {l}2}\log(-\lambda ).\tag5.12
$$ 

This ends the proof of the general assertion on the contributions from
the $S_i$.

The information can be sharpened further by considering the parity of
the terms in (5.3). Here the symbols of $L_m$ and $A^{2\nu }$ have
even-even alternating parity, whereas the
symbol of $\tfrac A {|A'|}$ has even-odd alternating parity (cf\. Definition
5.1ff.), so the symbol 
$b^{m,\nu }$ of the composition 
$B^{m,\nu }=L_mA^{2\nu }\tfrac A {|A'|}$ has even-odd alternating parity.
Since the integral over the sphere $\{|\xi '|=1\}$ of an odd function
vanishes, and $b^{m,\nu }_{1-n}$ is odd in $\xi '$ when $n$ is odd, we conclude that$$
c_{m+2\nu +n-1}(x')=0 \text{ if $n$ is odd.}\tag5.13
$$
So when $n$ is odd, there are no logarithmic contributions at all!
\qed
\enddemo 

This leads to:

\proclaim{Theorem 5.3} Let $D$ be a perturbation of $D^0$ as in
{\rm (1.1)} near $X'$ such that $P_0$ and all terms in the Taylor
expansion of $P_1$ in $x_n$ commute with $A$.  Let $F$ be a
differential operator in $E_1$ of order $m' $, and let $r+1>\frac{n+m' }2$.
If $n$ is odd,  
the resolvent and heat operator, resp\. zeta function, associated with $\Delta _B$ have
trace expansions without logarithms, resp\. meromorphic extensions
without double poles:
$$\aligned
\Tr (F \partial _ \lambda ^r(\Delta _B-\lambda )^{-1})&\sim
\sum_{-n\le k< \infty } \tilde a_{ k}(-\lambda )^{\frac{m' -k}2-r-1}
,\\
\Tr (F  e^{-t\Delta _B})&\sim
\sum_{-n\le k<\infty } a_{ k}t^{\frac{k-m' }2} 
,\\
\Gamma (s)\zeta (F,\Delta _B,s)\equiv \Gamma (s)\Tr (F \Delta _B^{-s})&\sim
\sum_{-n\le k<\infty } \frac{a_{ k}}{s+\frac{k-m'}2}-\frac{\Tr(F \Pi
_0(\Delta _B))}{s} 
,
\endaligned\tag 5.14$$
where the coefficients are locally determined for $-n\le k<0$ (for
$-n\le k<m'$ if $F$ is tangential. Here $\tilde a_{-n}$ and $a_{-n}$
vanish if $m'$ is odd.
\endproclaim 

\demo{Proof} Recall (1.21). We have from Theorem 3.3
that there is an expansion
$$
\Tr (F \partial _ \lambda ^r(\Delta _B-\lambda )^{-1})\sim
\sum_{-n\le k< 0 } \tilde c_{ k}(-\lambda )^{\frac{m' -k}2-r-1} 
+\sum_{k\ge 0}\bigl({ \tilde c'_{ k}}\log (-\lambda
)+{\tilde c''_{ k}}\bigr)(-\lambda )^{\frac 
{m' -k}2-r-1},\tag 5.15$$
where the $\tilde c_k$ and $\tilde c'_k$ are locally determined (and $\tilde c_{-n}$
vanishes if $m'$ is odd); we have to show that all the
$\tilde c'_k$ vanish.
For any $K$ we can expand $x_nP_1$ as in  (4.1), and since we know from
Theorem 3.5 that removing the remainder can only affect the logarithmic terms
with $k\ge K$, this reduces the problem to the case where $D$ is as
in (4.2). 
We apply Theorem 4.7. The $\psi $do $Q^r_{\lambda ,+}$ gives no
log-terms. The s.g.o\. $G'_{M,r}$ 
contributes with a trace expansion that is  $O(\ang\lambda ^{-M'})$,
where we can get $M'$ as large as we want by taking $M$ large.
Finally, $\Tr_X(\zeta _{\varepsilon  }G_{M,r}\zeta _{\varepsilon
 })=\Tr_{X^0}G_{M,r}+O(\ang\lambda ^{-N})$, any $N$, where $G_{M,r}$
is a finite sum of terms with structure as 
described in
(4.48)--(4.49). Let us decompose and  Taylor expand $F$:
$$
\aligned
F&=\sum_{m=0}^{m'}F_{m}(x,D_{x'})\partial
_{x_n}^m=F_{(l_0)}+x_n^{l_0+1}F'_{(l_0)},\text{ where}\\
F_{(l_0)}&=
\sum_{m=0}^{m'}\sum_{l=0}^{l_0}x_n^lF_{m,l}(x',D_{x'})\partial
_{x_n}^m,
\endaligned
$$
and $F'_{(l_0)}$ is of order $m'$. As usual,  the trace resulting from the
remainder $x_n^{l_0+1}F'_{(l_0)}$ has log-powers beginning at an index that goes to $\infty $
when $l_0\to\infty $. In the finite sum, the $x_n^l\partial _{x_n}^m$
applied to $x_n^kK_{A_\lambda }$ give other linear combinations of
terms $x_n^{k'}K_{A_\lambda }$, so when we take the normal trace of
$F_{(l_0)}G_{M,r}$, we get a finite sum of terms as in (4.49). By
Theorem 5.2, they contribute no logarithmic terms. Thus all the
$\tilde c'_k$ in (5.15) are zero. 

The nonlocal coefficients start in general at the
power $(-\lambda )^{\frac{m'}2-r-1}$, but when $F$ is tangential,
they start at the power $(-\lambda )^{-r-1}$ (cf\. Theorem 3.3, also
for the statement on $\tilde a_{-n}$).

This carries over 
to a heat trace expansion and a zeta function expansion by the
transitions explained e.g\. in [GS96], cf\. Corollary 3.7 above.
\qed 
\enddemo

In particular, this extends qualitatively the result of [GS96] on the
product case 
to trace expansions where the $x_n$-independent morphism $\varphi ^0$
considered there is replaced by an arbitrary differential
operator $F$.

Taking in particular  $F=F_1\psi D$, where $\psi $ is a morphism from
$E_2$ to $E_1$, and $F_1$ is a differential operator in $E_1$, we extend the results to eta-expansions (using that
\linebreak$(\Delta _B-\lambda )^{-r-1}$ maps into the domain of $D_\ge$, and
$ D$ vanishes on $V_0(\Delta _B)$):

\proclaim{Corollary 5.4} Assumptions on $D$ as in Theorem {\rm
5.3}. Let  $\psi $ be a morphism from $E_2$ to $E_1$ and $F_1$ a
differential operator in $E_1$ of order $m'$, and let
$r+1>\frac{n+m'+1 }2$. If $n$ is odd, there are expansions without
logarithms resp\. double poles:
$$\aligned
\Tr (F_1\psi D_\ge \partial _ \lambda ^r(\Delta _B-\lambda )^{-1})&\sim
\sum_{-n\le k< \infty } \tilde b_{ k}(-\lambda )^{\frac{m'+1 -k}2-r-1}
,\\
\Tr (F_1\psi D_\ge  e^{-t\Delta _B})&\sim
\sum_{-n\le k<\infty } b_{ k}t^{\frac{k-m'-1 }2} 
,\\
\Gamma (s)\eta (F_1\psi ,D_\ge,2s-1)\equiv \Gamma (s)\Tr (F_1\psi
D_\ge\Delta _B^{-s})&\sim 
\sum_{-n\le k<\infty } \frac{b_{ k}}{s+\frac{k-m'-1}2} 
.
\endaligned\tag 5.16$$
The coefficients $b_k$ are locally determined for $-n\le k<0$ (for
$-n\le k<m'$ if $F_1$ is tangential, cf\. Theorem {\rm 2.10}).
The coefficients $\tilde b_{-n}$ and $b_{-n}$ vanish if $m'$ is even.
\endproclaim

Let us also make some observations on the case where $n$ is even.
One can ask whether the expansions in general will have
logarithms at both integer and half-integer powers for $k>0$. As shown
in [GS96], this is
not so for the product case with a factor $\varphi ^0$, where the
terms with $k$ even $>0$ vanish. The result in [GS96] was based on explicit
trace
expansions of the zeta and eta functions of $A$; we can
now show by our qualitative arguments that the result extends
structurally to
tangential $x_n$-independent factors $F$.

\proclaim{Theorem 5.5} Let $D$ be of product type and let $F$ be a
differential operator of order $m' $. Let $r+1>\frac{n+m'}2$. When
$n$ is odd, the trace expansions are without 
logarithmic terms as in {\rm (5.14)}. When $n$ is even and $F$ is
tangential and $x_n$-independent near $X'$, then the trace expansions
have logarithms only at the ``zero'th'' power and at subsequent half-integer
powers
(resp\. double poles only at zero and at negative half-integers):
$$\align
\Tr (F \partial _ \lambda ^r(\Delta _B-\lambda )^{-1})&\sim
\sum_{-n\le k< \infty } \tilde a_{ k}(-\lambda )^{\frac{m'
-k}2-r-1}+{ \tilde c_0}(-\lambda )^{-r-1}\log (-\lambda
)\\
&\qquad
+\sum_{j\ge 0}{ \tilde c_{2j+1}}(-\lambda )^{-j-\frac12-r-1}\log (-\lambda
)
,\\
\Tr (F  e^{-t\Delta _B})&\sim
\sum_{-n\le k<\infty } a_{ k}t^{\frac{k-m' }2}+{ c_0}\log t
+\sum_{ j\ge 0}{ c_{2j+1}}t^{j+\frac 12}\log t 
,\tag 5.17\\
\Gamma (s)\Tr (F \Delta _B^{-s})&\sim
\sum_{-n\le k<\infty } \frac{a_{ k}}{s+\frac{k-m'}2}
+\frac{- c_0}{s^2}-\frac{\Tr(F \Pi _0(\Delta _B))}s
+\sum_{ j\ge 0}\frac{- c_{2j+1}}{(s+j+\frac 12)^2}.
\endalign
$$
The following coefficients are locally determined: $\tilde a_k$ and
$a_k$ for $-n\le k<m'$, $\tilde c_0$ and $c_0$, $\tilde c_{2j+1}$ and
$c_{2j+1}$ for $j\ge 0$.
\endproclaim 

\demo{Proof} The statement for $n$ odd has already been shown in
Theorem 5.3, so let
$n$ be even and let $F$ be tangential and $x_n$-independent near $X'$.
We have from Theorem 3.3 that there is an expansion (3.6) with $\mu
=(-\lambda )^{\frac12}$, so the point is now to show that only
certain log-terms appear. Since we are in the product case,
the normal trace of the s.g.o\. part  of $F\partial _\lambda
^rR^0_\lambda $ is (modulo contributions whose traces are $O(\ang{\lambda
}^{-N})$ for any $N$):
$$
\tr_n F\partial _\lambda ^r G^0_{\lambda }=
F\partial _\lambda ^r\bigl(\tfrac{-1}{4\lambda
}\bigl[\tfrac{A^2}{ A_\lambda 
^2}+\tfrac{A}{A_\lambda }\bigr] +\tfrac{A}{4\lambda
A_\lambda }\tfrac A{|A'|} 
+\tfrac A{4\lambda |A'| }+\tfrac1{4\lambda }\Pi _0\bigr),\tag 5.18
$$
as can be deduced from [GS96, (3.9)].
We appeal to 
Theorem 5.2 and its proof. In (5.18), the strongly polyhomogeneous terms in
[\dots] and the term with $\Pi _0$  
give no logs. The term$$
F\partial _\lambda ^r\tfrac{A}{4\lambda
A_\lambda }\tfrac A{|A'|} \tag5.19
$$
gives, when $\partial _\lambda ^r$ is carried out,  a linear
combination of terms (5.5) with $l-j$ odd (equal to 
$-3-2r$), hence it gives
log-power terms at half-integers, $c(x')(-\lambda )^{-\frac32-r-\nu }\log(-\lambda )$,
$\nu =0,1,2,\dots$, besides pure power terms. 
The term 
$$F\partial _\lambda ^r\frac {
A}{4\lambda |A'|}\tag5.20
$$ has the form with $j=0$ analyzed in the proof,
so it gives rise to 
one logarithmic term $c(x')(-\lambda )^{-r-1}\log(-\lambda )$ (as in (5.12))
besides pure power terms. 

This shows the first formula in (5.17), written with a different
enumeration convention than in (3.6) for the logarithmic terms. The
other formulas follow as in Corollary 3.7.
\qed
\enddemo 


This extends the result in [GS96]. One can analyze the eta function
in a similar way.

We recall from [G92], [GS96] that when $F=1$, the log-term at the
power $-r-1$ (resp\. at the power zero in the heat expansion)
vanishes. Note that the log-terms at the other powers stem
from one single log-producing term (5.19) that gives {\it odd} values
of $l-j$. If
$F$ is allowed to depend on $x_n$, the powers of $x_n$ in its Taylor
expansion will give rise to terms like (5.19) multiplied by negative
powers of $A_\lambda $, and terms like (5.20) multiplied with
negative powers of $A_\lambda $ (cf\. Lemma 4.3 (iii)). Then
both even and odd negative powers of $A_\lambda $ 
will occur, giving series of log-terms both with half-integer and with
integer powers. --- Also the effect of normal derivatives in $F$ can be
discussed in this way.

Now let us turn to non-product cases. Here, even when $F=1$, one can
expect many integer and half-integer log-terms,
as described in
the following remark.

\example {Remark 5.6} (Even $n$.) Consider the case of a nonzero
perturbation as 
in (4.2) and let just $F=1$. A thorough analysis seems 
unmanageable at this point, but we can get some evidence for what to
expect by analyzing the second term in the expansion of $\Cal R^0_M$ in
(4.35). By circular permutation,
$$
\Tr(\partial _\lambda ^r(\Cal R^0\zeta _\varepsilon \Cal P\Cal R^0))
=\Tr(\partial _\lambda ^r(\zeta _\varepsilon \Cal P(\Cal R^0)^2))=
-\Tr(\zeta _\varepsilon \Cal P\partial _\lambda ^{r+1}\Cal R^0).\tag5.21
$$
Consider the log-terms produced by its singular Green part. 

If
$\Cal P=\Cal P_0$, we can apply Theorem 5.5 directly to see that
there is an expansion with logs at half-integer powers $-l-\frac12-r-2$ and a
single log-term at the integer power $-r-2$.

If $\Cal P$ contains a term $x_n^k\Cal P_k$ with $k\ge 1$, the normal
trace of the  
s.g.o\. contribution 
from this term will be similar to (5.18), of the form:
$$
\tr_n Lx_n^k\partial _\lambda ^{r+1} G^0_{\lambda }=
k!L\partial _\lambda ^{r+1}\bigl(\tfrac1{(2A_\lambda
)^k}\bigl(\tfrac{-1}{4\lambda }\bigl[\tfrac{A^2}{ A_\lambda
^2}+\tfrac{A}{ A_\lambda }\bigr] +\tfrac{A}{4\lambda
A_\lambda }\tfrac A{|A'|} 
+\tfrac A{4\lambda |A'| }+\tfrac1{4\lambda }\Pi _0\bigr)\bigr),\tag 5.22
$$
in view of Lemma 4.3 (iii).
Again, the terms in $[\dots]$ produce no logs, but the interesting
fact is that now the terms generalizing (5.19) and (5.20) together
contain $A_\lambda $ in both even and odd negative powers
(and $\lambda $ in only integer powers), so that by Theorem 5.2,
log-terms are 
produced at both integer and half-integer powers from a certain step
on.
So already the second term in (4.35) will then
contribute sequences of nontrivial log-power terms both with integer
and half-integer powers; a strong 
indication that such a structure will be found in general.
\endexample

\enddocument